\documentclass[12pt]{amsart}
\begin{document}
\theoremstyle{plain}
\newtheorem{thm}{Theorem}[section]
\newtheorem*{thm*}{Theorem}
\newtheorem{prop}[thm]{Proposition}
\newtheorem*{prop*}{Proposition}
\newtheorem{lemma}[thm]{Lemma}
\newtheorem{cor}[thm]{Corollary}
\newtheorem*{conj*}{Conjecture}
\newtheorem*{cor*}{Corollary}
\newtheorem{defn}[thm]{Definition}
\theoremstyle{definition}
\newtheorem*{defn*}{Definition}
\newtheorem{rems}[thm]{Remarks}
\newtheorem*{rems*}{Remarks}
\newtheorem*{proof*}{Proof}
\newtheorem*{not*}{Notation}
\newcommand{\npartial}{\slash\!\!\!\partial}
\newcommand{\Heis}{\operatorname{Heis}}
\newcommand{\Solv}{\operatorname{Solv}}
\newcommand{\Spin}{\operatorname{Spin}}
\newcommand{\SO}{\operatorname{SO}}
\newcommand{\ind}{\operatorname{ind}}
\newcommand{\Index}{\operatorname{index}}
\newcommand{\ch}{\operatorname{ch}}
\newcommand{\rank}{\operatorname{rank}}
\newcommand{\abs}[1]{\lvert#1\rvert}
 \newcommand{\A}{{\mathcal A}}
        \newcommand{\D}{{\mathcal D}}\newcommand{\HH}{{\mathcal H}}
        \newcommand{\LL}{{\mathcal L}}
        \newcommand{\B}{{\mathcal B}}
        \newcommand{\K}{{\mathcal K}}
\newcommand{\oo}{{\mathcal O}}
         \newcommand{\PP}{{\mathcal P}}
        \newcommand{\s}{\sigma}
        \newcommand{\coker}{{\mbox coker}}
        \newcommand{\p}{\partial}
        \newcommand{\dd}{|\D|}
        \newcommand{\n}{\parallel}  
\newcommand{\bma}{\left(\begin{array}{cc}}
\newcommand{\ema}{\end{array}\right)}
\newcommand{\bca}{\left(\begin{array}{c}}
\newcommand{\eca}{\end{array}\right)}
\newcommand{\sr}{\stackrel}
\newcommand{\da}{\downarrow}
\newcommand{\tD}{\tilde{\D}}
        \newcommand{\R}{\mathbf R}
        \newcommand{\C}{\mathbf C}
        \newcommand{\h}{\mathbf H}
\newcommand{\Z}{\mathbf Z}
\newcommand{\N}{\mathbf N}
\newcommand{\tto}{\longrightarrow}
\newcommand{\ben}{\begin{displaymath}}
        \newcommand{\een}{\end{displaymath}}
\newcommand{\be}{\begin{equation}}
\newcommand{\ee}{\end{equation}}

        \newcommand{\bean}{\begin{eqnarray*}}
        \newcommand{\eean}{\end{eqnarray*}}
\newcommand{\nno}{\nonumber\\}
\newcommand{\bea}{\begin{eqnarray}}
        \newcommand{\eea}{\end{eqnarray}}
\newcommand{\supp}[1]{\operatorname{#1}}
\newcommand{\norm}[1]{\parallel\, #1\, \parallel}
\newcommand{\ip}[2]{\langle #1,#2\rangle}
\setlength{\parskip}{.3cm}
\newcommand{\nc}{\newcommand}
\nc{\nt}{\newtheorem}
\nc{\gf}[2]{\genfrac{}{}{0pt}{}{#1}{#2}}
\nc{\mb}[1]{{\mbox{$ #1 $}}}
\nc{\real}{{\mathbb R}}
\nc{\comp}{{\mathbb C}}
\nc{\ints}{{\mathbb Z}}
\nc{\Ltoo}{\mb{L^2({\mathbf H})}}
\nc{\rtoo}{\mb{{\mathbf R}^2}}
\nc{\slr}{{\mathbf {SL}}(2,\real)}
\nc{\slz}{{\mathbf {SL}}(2,\ints)}
\nc{\su}{{\mathbf {SU}}(1,1)}
\nc{\so}{{\mathbf {SO}}}
\nc{\hyp}{{\mathbb H}}
\nc{\disc}{{\mathbf D}}
\nc{\torus}{{\mathbb T}}
\newcommand{\tk}{\widetilde{K}}
\newcommand{\boe}{{\bf e}}\newcommand{\bt}{{\bf t}}
\newcommand{\vth}{\vartheta}
\newcommand{\CGh}{\widetilde{\CG}}
\newcommand{\db}{\overline{\partial}}
\newcommand{\tE}{\widetilde{E}}
\newcommand{\tr}{\mbox{tr}}
\newcommand{\ta}{\widetilde{\alpha}}
\newcommand{\tb}{\widetilde{\beta}}
\newcommand{\txi}{\widetilde{\xi}}
\newcommand{\hV}{\hat{V}}
\newcommand{\IC}{\mathbf{C}}
\newcommand{\IZ}{\mathbf{Z}}
\newcommand{\IP}{\mathbf{P}}
\newcommand{\IR}{\mathbf{R}}
\newcommand{\IH}{\mathbf{H}}
\newcommand{\IG}{\mathbf{G}}
\newcommand{\CC}{{\mathcal C}}
\newcommand{\CD}{{\mathcal D}}
\newcommand{\CS}{{\mathcal S}}
\newcommand{\CG}{{\mathcal G}}
\newcommand{\CL}{{\mathcal L}}
\newcommand{\CO}{{\mathcal O}}
\nc{\ca}{{\mathcal A}}
\nc{\cag}{{{\mathcal A}^\Gamma}}
\nc{\cg}{{\mathcal G}}
\nc{\chh}{{\mathcal H}}
\nc{\ck}{{\mathcal B}}
\nc{\cl}{{\mathcal L}}
\nc{\cm}{{\mathcal M}}
\nc{\cn}{{\mathcal N}}
\nc{\cs}{{\mathcal S}}
\nc{\cz}{{\mathcal Z}}
\nc{\sind}{\sigma{\rm -ind}}
\newcommand{\la}{\langle}
\newcommand{\ra}{\rangle}

\begin{center}
 \title{The local index formula in semifinite von Neumann algebras I: Spectral flow}

 \vspace{.5 in}

\author{}
\maketitle
{\bf Alan L. Carey}\\Mathematical Sciences Institute\\
Australian National University\\
Canberra, ACT. 0200, AUSTRALIA\\
e-mail: acarey@maths.anu.edu.au\\
\vspace{.2 in}

{\bf John Phillips}\\Department of Mathematics and Statistics\\
University of Victoria\\Victoria, B.C. V8W 3P4, CANADA\footnote{Address for
correspondence}\\
e-mail: phillips@math.uvic.ca\\
\vspace{.2 in}

{\bf Adam Rennie}\\
School of Mathematical and Physical Sciences\\
University of Newcastle\\
Callaghan, NSW, 2308 AUSTRALIA\\
e-mail: adam.rennie@newcastle.edu.au\\

\vspace{.2 in}

{\bf Fyodor A. Sukochev}\\
School of Informatics and Engineering\\
Flinders University\\
Bedford Park S.A 5042 AUSTRALIA\\
e-mail: sukochev@infoeng.flinders.edu.au\\
\vspace{.25 in}

All authors were supported by grants from ARC (Australia) and
NSERC (Canada), in addition the third named author acknowledges a 
University of Newcastle early career researcher grant and the first
named author acknowledges the Clay Mathematics Institute for whom this 
research was begun and support from the Erwin Schr\"{o}dinger Institute's
Noncommutative Geometry program.

\end{center}

\newpage
\centerline{{\bf Abstract}}
We generalise the local index formula of Connes and Moscovici to the
case of spectral triples for a $*$-subalgebra $\A$
of a general semifinite von Neumann algebra. In this setting it 
gives a formula for spectral flow along a path joining an unbounded self 
adjoint Breuer-Fredholm operator, affiliated to the von Neumann
algebra, to a unitarily equivalent operator.
Our proof is novel even in the setting of the original theorem 
and relies on the introduction of a function valued cocycle
which is `almost'
a $(b,B)$-cocycle in the cyclic cohomology of $\mathcal A$.
\footnote{
von Neumann algebra, Fredholm module, cyclic cohomology, chern character,
spectral flow.}

\newpage       
\allowdisplaybreaks
\section{Introduction}

The odd local index theorem of Connes and Moscovici \cite{CM}
may be thought of as a far reaching generalisation of the classical index 
theorem for Toeplitz operators. It is thus a natural question 
to ask whether the index theorem of
Coburn, Douglas, Schaeffer and Singer \cite{CDSS, CMX} proved in the 
setting of semifinite von Neumann algebras and giving a
topological formula for the 
 Breuer-Fredholm index of Wiener-Hopf operators
with almost periodic symbol is the prototype for a von Neumann algebra version 
of the local index theorem.
This question was answered in the affirmative by our noncommutative geometry
calculation of the index of Toeplitz operators with
noncommutative symbol \cite{CPS2, L, PR}. 
In both cases there is
a clear interpretation of the index as computing spectral flow along a
certain path of unbounded self-adjoint Breuer-Fredholm operators.
(This follows from recent work of some of us in \cite{CP1, CP2}
interpreting the Breuer-Fredholm index of the \cite{CDSS}
Wiener-Hopf operators as `type II' spectral flow.)
This motivated the present general study of the local index 
formula of Connes and Moscovici
in the setting of semifinite von Neumann algebras
via a computation of spectral flow 
along a path of self-adjoint unbounded Breuer-Fredholm operators.

This line of reasoning touches on a more general program outlined
by \cite{BeF} for developing a theory
of  `von Neumann spectral triples'. 
In addition to \cite{CDSS}, examples which suggest that 
this has interest include differential
operators with almost periodic symbol \cite{Sh},
the $L^2$-index theorem (see \cite{M} and references therein),
foliations \cite{Co4, BeF, Pr} as well as the example
of spectral flow cited above \cite{CP1,CP2}.

The starting point is a Hilbert space $\mathcal H$ on which there is
an unbounded densely defined self adjoint
Breuer-Fredholm operator $\mathcal D$.
In this setting
Carey-Phillips introduced an integral formula for the spectral flow 
along the linear path joining 
$\D$ to a unitarily equivalent operator $u\D u^*$, \cite{CP1,CP2}, which we 
present in Equation (\ref{basicformula}) in Section \ref{spectralflow}. 
The natural framework for this formula is that of
 odd spectral triples (generalised
to the von Neumann setting as  in \cite{CP1,BeF,CPS1,CPS2})
so that $u$ is an element of a $*$-subalgebra $\mathcal A$ of
a semifinite von Neumann algebra $\mathcal N$ acting 
non-degenerately on $\mathcal H$.

Given the analytic formula of \cite{CP1,CP2} 
for the spectral flow, our task is to show that there is an 
equality with
a cohomological formula. We derive the cohomological formula in several steps. 

The first step, described in Section \ref{G}, exploits the fact 
that the analytic formula gives spectral flow as the integral of an exact 
one-form, \cite{CP1,CP2}. 
The exactness of the one-form allows us to change the path of integration 
to obtain a new formula which is amenable to perturbation theory methods. 
In Section \ref{secresolvent} we employ a perturbation expansion of the 
resolvent to write spectral flow in terms of a `function-valued cochain' in 
the $(b,B)$ bicomplex of cyclic cohomology. Our function-valued cochain  is 
reminiscent of, though distinct from, Higson's `improper cocycle' \cite{hig}. 
Our cochain is a cocycle modulo functions holomorphic in a half-plane. 
We refer to this cochain as the resolvent cocycle and it should be thought of as
a substitute for the JLO cocycle (which is the starting point for the 
argument of Connes-Moscovici).

The resolvent cocycle can be further expanded employing the quantised 
pseudodifferential calculus of Connes-Moscovici, \cite{CM}. This is done 
in Section \ref{residue} using material developed in Section \ref{psido}.
The end result is an expression for the 
spectral flow in terms of a sum of generalised zeta functions of the form
\ben \zeta_b(z)=\tau(b(1+\D^2)^{-z}),\ \ \ b\in\mathcal N.\een
This sum of zeta functions is meromorphic in a half-plane, with (at worst) 
only a single simple pole in this half-plane. 
The residue at this pole is precisely the spectral flow. 

Under the assumption that the individual zeta functions in this sum
analytically continue to a deleted neighbourhood of the critical point,
we may take residues of the individual terms at the critical point. 
The resulting formula, when $\mathcal N={\mathcal B}({\HH})$,
 is essentially that which is obtained by pairing the (odd, 
renormalised) 
cyclic cocycle obtained by Connes and Moscovici, \cite{C4,CM}, with 
the Chern character of the unitary $u^*$. The only difference between the 
two formulae is that we do not assume $\D$ is invertible and hence use inverse 
powers of $(1+\D^2)^{1/2}$, whereas Connes-Moscovici assume that $\D$ is 
invertible and use inverse powers of $\dd$. 

The novel aspects of our approach are:

$\bullet$ Our result calculates spectral flow in
semifinite von Neumann spectral triples generalising
part of the type I theory of \cite{CM}. 
Specifically, our formula for spectral 
flow is given in terms of a cyclic cocycle, 
which is the generalisation to semi-finite von Neumann algebras of
the residue cocycle of \cite{CM}. 
This provides an extension of \cite{CPS2, L, PR}.

$\bullet$ Only the final step of our proof requires the analytic continuation 
property of the generalised zeta functions. Indeed, we express spectral 
flow as the {\em residue} of a sum of zeta functions without invoking {\em 
any} analytic continuation hypothesis. 

$\bullet$ Assuming the individual zeta 
functions in the above sum have analytic continuations to a deleted 
neighbourhood of the 
critical point allows us to write spectral flow as a sum of residues of 
zeta functions. The residues of these zeta functions then assemble to form a 
$(b,B)$ cocycle for the algebra $\A$ of the spectral triple.

For examples with `dimension' less than $2$, these last two statements are true 
without {\em any} analytic continuation property.

$\bullet$ We make no assumptions on the decay of our zeta functions 
along vertical lines in the complex plane thus reducing the side conditions 
that need to be checked when applying the local index formula of \cite{CM}.

$\bullet$ Our proof
that the residue cocycle \cite{CM} is indeed a 
$(b,B)$-cocycle is quite simple even in the
general semifinite case by virtue of using our resolvent cocycle.

$\bullet$ Except for the need to verify a number of estimates, the strategy of 
our proof is straightforward, and applies to both the type I and type II cases.

$\bullet$ The idea of this proof in the odd case
can be adapted to handle the even case of 
the local index formula. The starting point for the even case is a generalised
McKean-Singer formula and the argument is presented
in  Part II.

$\bullet$ We remark that there is an unrenormalised version of the residue 
cocycle
in \cite{CM} containing an 
infinite number of terms 
in the case that one of the terms in the expansion has 
an essential singularity, whereas their renormalised version always has a 
bounded number of terms. The unrenormalised version presents an issue of 
convergence which is difficult to address. 
Since we do not pass through an intermediate step where the 
cocycle contains a potentially infinite number of terms, we are free to allow 
essential singularities from the outset.

In this paper we do not address the relationship of the residue cocycle to the 
Chern
character as is done in \cite{CM} leaving that to another place as it is not a
trivial step. Our residue cocycle necessarily 
involves zeta functions of $(1+D^2)^{-1/2}$ 
because $|D|$ may have zero in its continuous spectrum. This means 
the usual transgression arguments do not immediately work. 
It also means that the formula obtained when $\mathcal N={\mathcal B}({\HH})$
is a modification of a consequence
(Corollary II of \cite{CM}) of the `renormalised' version of the
Connes-Moscovici Local Index Theorem.

The  plan of the paper is that we group preliminary material
including notation and definitions in Section 2. Section 3 describes integral 
formulae for Spectral Flow, and the relation to cyclic cohomology.
The statement of the main result is in Section 4, and at this point we 
also provide an outline of the proof for the reader's convenience.

Section 5 establishes
the key estimates needed for the rest of the proof.
Section 6 summarises what we
need to know about the quantised pseudodifferential calculus.
In Section 7 
we derive our resolvent cocycle. These earlier ingredients come together in 
Section 8 to prove the 
main theorem.

\noindent{\bf Acknowledgements}
We thank Nigel Higson for discussions on his approach to the local index 
theorem and in particular on his point of view of the pseudodifferential
calculus.

\section{Definitions and Background}

\subsection{Semifinite Spectral Triples} 

We have adopted a notational convention correlated to the context in which
we are working.
A calligraphic $\D$ will always denote 
an unbounded self-adjoint operator forming part of a 
semifinite spectral triple $(\A,\HH,\D)$. A roman $D$ will 
denote a self-adjoint operator on a Hilbert space, usually with some side 
conditions. Later, starting from a semifinite
spectral triple $(\A,\HH,\D)$ in section \ref{G}, we construct
from the space $\HH$ on which 
a semifinite von Neumann algebra $\mathcal N$ acts, a 
new Hilbert space $\tilde{\HH}$, an algebra $\tilde{\mathcal N}$ 
and an operator $\tD$ on $\tilde{\HH}$ affiliated with 
$\tilde{\mathcal N}$. 
While this is an odd spectral triple (for an algebra containing $\mathcal A$), 
for us it is just a computational device inspired by ideas of \cite{G}. 
We begin with some semifinite versions of standard definitions and results.
Let ${\mathcal K}_{\mathcal N }$ be the 
$\tau$-compact operators in ${\mathcal N}$
(that is the norm closed ideal generated by the projections
$E\in\mathcal N$ with $\tau(E)<\infty$). Here $\tau$ is a fixed faithful 
semifinite trace on the von Neumann algebra ${\mathcal N}$.

\begin{defn} A semifinite
spectral triple $(\A,\HH,\D)$ is given by a Hilbert space $\HH$, a 
$*$-algebra $\A\subset \cn$ where $\cn$ is a semifinite von Neumann algebra 
acting on
$\HH$, and a densely defined unbounded self-adjoint operator $\D$ affiliated 
to $\cn$ such that

1) $[\D,a]$ is densely defined and extends to a bounded operator for all 
$a\in\A$

2) $(\lambda-\D)^{-1}\in\K_\cn$ for all $\lambda\not\in{\R}$
\end{defn}

\noindent{\bf Note}. In this paper, for simplicity of exposition,
we will deal only with unital algebras
${\mathcal A}\subset \cn$ where the identity of $\mathcal A$ is that of $\cn$.
Henceforth we omit the term semifinite
as it is implied by the use of a faithful normal semifinite trace $\tau$
on $\cn$ in all of the subsequent text. In this paper (Part I) we will only 
deal with odd spectral 
triples, \cite{Co4}, since spectral flow 
is the pairing of $K$-homology with $K_1$.
 
\begin{defn}\label{qck} A semifinite spectral triple $(\A,\HH,\D)$ is $QC^k$ 
for $k\geq 1$ 
($Q$ for quantum) if for all $a\in\A$ 
the operators $a$ and $[\D,a]$ are in the domain of $\delta^k$, where 
$\delta(T)=[\dd,T]$ is the partial derivation on $\cn$ defined by $\dd$. We 
say that 
$(\A,\HH,\D)$ is $QC^\infty$ if it is $QC^k$ for all $k\geq 1$.
\end{defn}

{\bf Note}. The notation is meant to be analogous to the classical case, but 
we introduce 
the $Q$ so that there is no confusion between quantum differentiability of 
$a\in\A$ and classical differentiability of functions.

\noindent{\bf Remarks concerning derivations, commutators and topology}.  By 
partial 
derivation we mean that $\delta$ is defined on some subalgebra of $\cn$ which 
need not be (weakly) dense in $\cn$. More precisely, $\mbox{dom}\delta=
\{T\in\cn:\delta(T)\mbox{ is bounded}\}$. We 
also note that if $T\in{\mathcal N}$, one can show that $[\dd,T]$ is bounded 
if and only if $[(1+\D^2)^{1/2},T]$ is bounded, by using the functional 
calculus to show that $\dd-(1+\D^2)^{1/2}$ extends to a bounded operator in 
$\cn$. In fact, writing $\dd_1=(1+\D^2)^{1/2}$ and $\delta_1(T)=[\dd_1,T]$ we 
have $dom(\delta^n)=dom(\delta_1^n)$ for all $n.$

\begin{proof} Let $f(\D)=(1+\D^2)^{1/2}-\dd$, so, as noted above, $f(\D)$ 
extends to a bounded operator in $\cn$. Since
$$ \delta_1(T)-\delta(T)=[f(\D),T]$$ is always bounded, $\mbox{dom}\delta=
\mbox{dom}\delta_1$. Now $\delta\delta_1=\delta_1\delta$, so
\begin{eqnarray*}
\delta^2_1(T)-\delta^2(T)&=&\delta_1(\delta_1(T))-\delta_1(\delta(T))+
\delta_1(\delta(T))-\delta(\delta(T))\\
&=&[f(\D),\delta_1(T)]+[f(\D),\delta(T)].\end{eqnarray*}
Both terms on the right hand side are bounded, so $\mbox{dom}\delta^2=
\mbox{dom}\delta_1^2$. The proof proceeds by induction.
\end{proof}

Thus the condition defining $QC^\infty$ can be replaced by
\ben a,[\D,a]\in\bigcap_{n\geq 0}\mbox{dom}\delta_1^n\ \ \ \forall a\in\A.\een
This is important as we {\em do not assume at any point that $\dd$ is 
invertible}. 

If $(\A,\HH,\D)$ is a $QC^\infty$ spectral triple, we may endow the algebra 
$\A$ with the topology determined by the seminorms
\ben a\tto \n \delta^k(a)\n+\n\delta^k([\D,a])\n,\ \ \ k=0,1,2,...\een
We call this topology the $\delta$-topology and observe that by 
\cite[Lemma 16]{R} we may, without loss of generality, suppose that $\A$ is 
complete in the $\delta$-topology by completing if necessary. This completion 
is Fr\'{e}chet and stable under the holomorphic functional calculus, so we 
have a sensible spectral theory and 
$K_*(\A)\cong K_*(\bar\A)$ via inclusion, where $\bar\A$ 
is the $C^*$-completion of $\A$.

Next we observe that if $T\in\cn$ and $[\D,T]$ is bounded, then $[\D,T]\in\cn$.

\begin{proof} Observe that $\D$ is affiliated with $\cn$, and so commutes with 
all projections in the commutant of $\cn$, and the commutant of $\cn$ preserves 
the domain of $\D$. Thus if $[\D,T]$ is bounded, it too commutes with all 
projections in the commutant of $\cn$, and these projections preserve the 
domain of $\D$, and so $[\D,T]\in\cn$. 
\end{proof}
Similar comments apply to $[\dd,T]$, $[(1+\D^2)^{1/2},T]$ and the more exotic 
combinations such as $[\D^2,T](1+\D^2)^{-1/2}$ which we will encounter later. 

Recall from \cite{FK} that if $S\in\mathcal N$, the {\bf t-th generalized
singular value} of S for each real $t>0$ is given by
$$\mu_t(S)=\inf\{||SE||\ \vert \ E \mbox{ is a projection in }
{\mathcal N} \mbox { with } \tau(1-E)\leq t\}.$$
The ideal $\LL^1({\mathcal N})$ consists of those 
operators $T\in {\mathcal N}$ such that
$ \n T\n_1:=\tau( |T|)<\infty$
where $|T|=\sqrt{T^*T}$.
In the Type I setting this is the usual trace class ideal. We will 
simply write $\LL^1$ for this ideal in order to simplify the notation, and 
denote the norm on $\LL^1$ by $\n\cdot\n_1$. An alternative definition in terms 
of 
singular values is that $T\in\LL^1$ if
$ \|T\|_1:=\int_0^\infty \mu_t(T) dt <\infty.$

Note that in the case where 
${\mathcal N}\neq{\mathcal B}({\mathcal H})$, 
$\LL^1$ need not be 
complete in this norm but it is complete in the norm $||.||_1 + ||.||_\infty$.
(where $||.||_\infty$ is the uniform norm).

\subsection{Dimension Spectrum}

Many naturally occurring spectral triples satisfy a summability condition,
and such conditions allow one to define interesting cocycles.
In particular, one can obtain representatives of the Chern character.
The (finite)  summability conditions give a half-plane where the function
\be z\mapsto \tau((1+\D^2)^{-z})\label{zeta}\ee
is well-defined and holomorphic. In \cite{C4,CM}, a stronger condition was
imposed in order to prove the Local Index Theorem. This condition not only
specifies a half-plane where the function in (\ref{zeta}) is holomorphic, 
but also that this function analytically continues to ${\C}$ minus some 
discrete set. We clarify this in the following definitions.

\begin{defn} Let $(\A,\HH,\D)$ be a $QC^\infty$ spectral triple. The
algebra $\B(\A)\subseteq \mathcal N$ is the algebra of polynomials generated by
 $\delta^n(a)$ and $\delta^n([\D,a])$ for $a\in\A$ and $n\geq 0.$
 A $QC^\infty$ spectral triple $(\A,\HH,\D)$ has {\bf discrete
dimension spectrum} $Sd\subseteq {\C}$ if $Sd$ is a discrete set and for all
$b\in\B(\A)$ the function $\tau(b(1+\D^2)^{-z})$
is defined and holomorphic for $Re(z)$ large, and analytically continues
to ${\C}\setminus Sd$. We say the dimension spectrum is
{\bf simple} if this zeta function has poles of order at most one for all
$b\in\B(\A)$, {\bf finite} if there is a $k\in{\N}$ such that the function 
has poles of order at most $k$ for all $b\in\B(\A)$
and {\bf infinite}, if it is not finite.
\end{defn}

Connes and Moscovici impose the discrete dimension spectrum assumption to
prove their version of the local index theorem.
In this paper we employ a weaker condition
(explained in the next section) that is implied
by the  discrete dimension spectrum assumption.

\section{Spectral Flow}\label{spectralflow}

To place our results in their proper setting we need some
background from \cite{Ph,Ph1,PR}. Let
$\pi:{\mathcal N}\to {\mathcal N}/{\mathcal K_{\cn}}$
be the canonical mapping. A Breuer-Fredholm operator is one that
maps to an invertible operator under $\pi$, \cite{PR}.
In the Appendix to \cite{PR}, the theory of Breuer-Fredholm operators  
for the case where $\cn$ is not a factor is developed in analogy with the 
factor case of Breuer, \cite{B1,B2}. In Part II of this work,
we develop this theory even further.
  As usual $D$ is an unbounded densely defined self-adjoint
Breuer-Fredholm operator  on $\HH$ (meaning $D(1+D^2)^{-1/2}$ is 
bounded and Breuer-Fredholm in $\mathcal N$)
with $(1+D^2)^{-1/2}\in {\mathcal K}_\mathcal N$.
For a unitary $u\in \mathcal N$ such that $[D,u]$ is a bounded operator,
 the path
$$D_t^u:=(1-t)\,D+tuDu^*$$
of unbounded self-adjoint
Breuer-Fredholm operators is continuous in the sense that
$$F_t^u:=D_t^u\left(1+(D_t^u)^2\right)^{-\frac{1}{2}}$$
is a norm continuous path of self-adjoint Breuer-Fredholm operators in 
$\mathcal N$ \cite{CP1}.
Recall that the Breuer-Fredholm index of
a Breuer-Fredholm operator $F$ is defined by
$$ind(F)=\tau(Q_{ker F})-\tau(Q_{coker F})$$
where $Q_{ker F}$ and $Q_{coker F}$ are the projections onto
the kernel and cokernel of $F$.

\begin{defn*}
If $\{F_t\}$ is a continuous path of self-adjoint Breuer-Fredholm
operators
in $\mathcal N$, then the definition of the {\bf spectral flow} of the
path, $sf(\{F_t\})$ is based on the following sequence of observations
in \cite{P1}:

\noindent 1. The function $t\mapsto sign(F_t)$ is typically discontinuous as is
the
projection-valued mapping $t\mapsto P_t=\frac{1}{2}(sign(F_t)+1).$

\noindent 2. However,  $t\mapsto \pi(P_t)$ is continuous.

\noindent 3. If $P$ and $Q$ are projections in $\mathcal N$ and
$||\pi(P)-\pi(Q)||<1$
then $PQ:Q\HH \to P\HH$ is a Breuer-Fredholm operator and so
$ind(PQ)\in {\R}$ is well-defined. This requires
\S 3 of Part II.

\noindent 4. If we partition the parameter interval of $\{F_t\}$ so
that the $\pi(P_t)$ do not vary much in norm on each subinterval
of the partition then $$sf(\{F_t\}):=\sum_{i=1}^n
ind(P_{t_{i-1}}P_{t_i})$$
is a well-defined and (path-) homotopy-invariant number which
agrees with the usual notion of spectral flow in the type $I_\infty$
case.

\noindent 5. For $D$ and $u$ as above, we define the {\it spectral flow} of 
the path $D_t^u:=(1-t)\,D+tuDu^*$ to be the spectral flow of the path 
$F_t$ where $F_t=D_t^u\left(1+(D_t^u)^2\right)^{-\frac{1}{2}}$. We denote this 
by
\ben sf(D,uDu^*)=sf(\{F_t\}),\een
and observe that this is an integer in the
${\mathcal N}={\mathcal B}({\mathcal H})$ case and a real
number in the general semifinite case.
\end{defn*}

Special cases of
spectral flow in a semifinite von Neumann algebra were discussed in
\cite{M,P1,P2}. 

Let $P$ denote the projection onto the nonnegative spectral subspace of
$D$.
The spectral flow along $\{D_t^u\}$ is equal to $sf(\{F_t\})$
and by \cite{CP1} this is the Breuer-Fredholm index of
$PuPu^*$. (Note that sign$F^u_1=2uPu^*-1$ and that
for this special path we have $P-uPu^*$ is compact so
$PuPu^*$ is certainly Breuer-Fredholm from 
$uPu^*\HH \to P\HH$.) 
Now, \cite[Appendix B]{PR}, we have ind$(PuPu^*)=$ ind$(PuP)$.

\subsection{Spectral Flow Formulae}

We now introduce the spectral flow formula of Carey and Phillips, 
\cite{CP1,CP2} from which our other formulae follow. 
This formula starts with a semifinite spectral triple $(\A,\HH,\D)$
and computes the spectral flow from $\D$ to $u\D u^*$, 
where $u\in\A$ is unitary with $[\D,u]$ bounded, in the case
where $(\A,\HH,\D)$ is of dimension $p\geq 1$.
Thus for any $n>p$ we have by Theorem 9.3 of \cite{CP2}:
\be sf(\D,u\D 
u^*)=\frac{1}{C_{n/2}}\int_0^1\tau(u[\D,u^*](1+(\D+tu[\D,u^*])^2)^{-n/2})dt,
\label{basicformula}\ee
with $C_{n/2}=\int_{-\infty}^\infty(1+x^2)^{-n/2}dx$. 
This real number $sf(\D,u\D u^*)$
recovers the pairing of the $K$-homology class
$[\D]$ of $\mathcal A$ with the $K_1({\mathcal A})$ class $[u]$
(see below).There is a geometric way to view this formula. It is shown in 
\cite {CP2} that
the functional $X\mapsto \tau(X(1+(\D+X)^2)^{-n/2})$ on ${\mathcal N}_{sa}$
determines an exact one-form
on an affine space modelled on  ${\mathcal N}_{sa}$. Thus (\ref{basicformula})
represents the integral of this one-form along the path 
$\{\D_t=(1-t)\D+ tu\D u^*\}$
provided one appreciates that $\dot\D_t=u[\D,u^*]$ is a tangent vector
to this path.
Moreover this formula is 
scale invariant. By this we mean that if we replace $\D$ by $\epsilon\D$, for 
$\epsilon>0$, in the right hand side of (\ref{basicformula}), then the left 
hand 
side is unchanged, since spectral flow is invariant with respect to change of 
scale. We will use this fact later on at several points. 

\subsection{Relation to Cyclic Cohomology}

One can also interpret spectral flow (in the type I case) 
as the pairing between an odd $K$-theory 
class represented by a unitary $u$, and an odd $K$-homology class represented 
by 
$(\A,\HH,\D)$, \cite[Chapter III,IV]{Co4}. This point of view also makes sense 
in the 
general semifinite setting, though one must suitably interpret $K$-homology, 
\cite{CPRS1,CP2}.
A central feature of \cite{Co4} 
is the translation of the $K$-theory pairing to 
cyclic theory in order to obtain index theorems. One associates to a suitable 
representative 
of a $K$-theory class, respectively a $K$-homology class, a class in periodic 
cyclic homology, 
respectively a class in periodic cyclic cohomology, 
called a Chern character in both cases. The principal result is then  
\be sf(\D,u\D u^*)=\langle [u],[(\A,\HH,\D)]\rangle=-\frac{1}{\sqrt{2\pi i}}
\langle [Ch_*(u)],[Ch^*(\A,\HH,\D)]\rangle,\label{indpair}\ee
where $[u]\in K_1(\A)$ is a $K$-theory class with representative $u$ and 
$[(\A,\HH,\D)]$ is the $K$-homology class of the spectral triple $(\A,\HH,\D)$.
 
On the right hand side, $Ch_*(u)$ is the Chern character of $u$, and 
$[Ch_*(u)]$
its periodic cyclic 
homology class. Similarly $[Ch^*(\A,\HH,\D)]$ is the periodic cyclic 
cohomology class of 
the 
Chern 
character of $(\A,\HH,\D)$. {\em The analogue of Equation (\ref{indpair}), for 
a suitable cocycle associated to $(\A,\HH,\D)$, in the general semifinite case 
is part of our main result.}

We will not discuss periodic cyclic cohomology Chern characters in great detail
here, 
leaving that to another place. We do, however, require some of the basic 
definitions and results of cyclic theory, as well as one result on the Chern 
character of a unitary. We will use the normalised $(b,B)$-bicomplex
(see \cite{Co4,Lo}). 

We introduce the following 
linear spaces.
Let $C_m=\A\otimes \bar\A^{\otimes m}$
where $\bar\A$ is the quotient $\A/\C I$ with $I$ being the 
identity
element of $\A$ and (assuming with no loss of generality that $\A$ is complete 
in the $\delta$-topology) we employ the projective tensor product.
Let $C^m=Hom(C_m,\C)$ be the linear space of continuous multilinear functionals
on $C_m$.
We may define the $(b,B)$ bicomplex using these spaces 
(as opposed to $C_m=\A^{\otimes m+1}$ et cetera) and the resulting cohomology 
will be the same. This follows because the bicomplex defined using
$\A\otimes \bar\A^{\otimes m}$ is quasi-isomorphic to that defined 
using  $\A\otimes \A^{\otimes m}$.

A normalised
$\mathbf{(b,B)}${\bf-cochain}, $\phi$ is a finite collection of continuous 
multilinear 
 functionals on $\A$,
$$\phi=\{\phi_m\}_{m=1,2,...,M}\mbox{ with }\phi_m\in C^m.$$
It is a (normalised) $\mathbf{(b,B)}${\bf-cocycle} if, for all $m$,
$b\phi_m+B\phi_{m+2}=0$ where 
$b: C^m\to C^{m+1}$, $B:C^m\to C^{m-1}$
are the coboundary operators given by
$$(B\phi_m)(a_0,a_1,\ldots,a_{m-1})
=\sum_{j=0}^{m-1} (-1)^{(m-1)j}\phi_m(1,a_j,a_{j+1},\ldots,a_{m-1},a_0,
\ldots,a_{j-1})$$
$$(b\phi_{m-2})(a_0,a_1,\ldots,a_{m-1})=\hfill$$
$$\sum_{j=0}^{m-2}(-1)^j\phi_{m-2}(a_0,a_1,\ldots,a_ja_{j+1},\ldots,a_{m-1})
+(-1)^{m-1}\phi_{m-2}(a_{m-1}a_0,a_1,\ldots,a_{m-2})$$
 We write $(b+B)\phi=0$ for brevity.
Thought of as functionals on $\A^{\otimes m+1}$
a normalised cocycle will satisfy
$\phi(a_0,a_1,\ldots,a_n)=0$ whenever any $a_j=1$ for $j\geq 1$.
An {\bf odd} ({\bf even}) cochain has $\{\phi_m\}=0$ for $m$ even (odd).

 Similarly, a $\mathbf{(b^T,B^T)}${\bf-chain}, $c$ is a (possibly infinite) 
 collection 
$c=\{c_m\}_{m=1,2,...}$ with 
$c_m\in C_m$. 
The $(b,B)$-chain $\{c_m\}$ 
is a $\mathbf{(b^T,B^T)}${\bf-cycle} if $b^Tc_{m+2}+B^Tc_m=0$ for all $m$. More 
briefly, we write $(b^T+B^T)c=0$.
Here $b^T,B^T$ 
are the boundary operators of cyclic homology, 
and are the transpose of the coboundary operators 
$b,B$ in the following sense. 

The pairing between a $(b,B)$-cochain $\phi=\{\phi_m\}^M_{m=1}$ and a 
$(b^T,B^T)$-chain 
$c=\{c_m\}$ 
is given by
\ben \langle \phi,c\rangle= \sum_{m=1}^M\phi_m(c_m).\een
This pairing satisfies
\ben \langle (b+B)\phi,c\rangle=\langle\phi,(b^T+B^T)c\rangle.\een
We use this fact in Section 8 in the following way. 
We call $c=(c_m)_{m\ odd}$ an odd normalised
$\mathbf{(b^T,B^T)}${\bf-boundary} if 
there is some even chain $e=\{e_m\}_{m\ even}$ 
with $c_{m}=b^Te_{m+1}+B^Te_{m-1}$ for all $m$. 
If we pair a normalised $(b,B)$-cocycle $\phi$
with a normalised $(b^T,B^T)$-boundary $c$ we find
\ben\langle \phi,c\rangle=
\langle \phi,(b^T+B^T)e\rangle=\langle (b+B)\phi,e\rangle=0.\een
There is an analogous definition in the case of even chains
 $c=(c_m)_{m\ even}$.
All of the cocycles we consider in this paper are in fact defined
as functionals on $\oplus_m\A\otimes \bar\A^{\otimes m}$.
Henceforth we will drop the superscript on  $b^T,B^T$ and just write $b,B$
for both boundary and coboundary operators
as the meaning will be clear from the context.

We recall that
the Chern character $Ch_*(u)$ of a unitary $u\in\A$ is the following (infinite)
collection of odd chains 
$Ch_{2j+1}(u)$ satisfying $bCh_{2j+3}(u)+BCh_{2j+1}(u)=0$,
\ben Ch_{2j+1}(u)=(-1)^jj!u^*\otimes u\otimes u^*\otimes\cdots\otimes u\ \ \ 
(2j+2\ \ \mbox{entries}) .\een
It is well known to experts that 
\be Ch_{*}(u^*)+Ch_{*}(u)\label{cob}\ee
is homologous to zero 
in the normalised (entire) $(b,B)$ chain complex
however an accessible argument eluded us so we present one 
here. A similar statement holds in the periodic theory.

\begin{lemma}\label{unitary}
For $u$ unitary in $\mathcal A$,
$Ch_{*}(u^*)+Ch_{*}(u)$ 
is a boundary in the odd normalised entire cyclic homology of $\mathcal A$. 
\end{lemma}

\begin{proof} That $Ch_*(u)$ defines an entire cycle is presented in \cite{G}.
We work with elements of $\A\otimes\bar{\mathcal A}^{\otimes n}$ 
written as $n+1$-tuples in the normalised version of the $(b,B)$ complex. 
So in the normalised complex with $n=2m+1$, $m=1,2,\ldots$
let 
$$w_{2m+1}=(u^{-1},u,\ldots,u^{-1},u)\in 
\A\otimes\bar{\mathcal A}^{\otimes 2m+1}$$ 
and then
$$Bw_{2m+1}=B(u^{-1},u,\ldots,u^{-1},u)
=(m+1)(1,u^{-1},u,\ldots,u^{-1},u)-(m+1)
(1,u,u^{-1}\ldots,u,u^{-1})$$
and $bw_{2m+1}\in  \A\otimes\bar{\mathcal A}^{\otimes 2m}$
$$bw_{2m+1}=b(u^{-1},u,\ldots,u^{-1},u)=
(1,u^{-1},u,\ldots,u^{-1},u)-(1,u,u^{-1},\ldots,u,u^{-1}).$$
Thus $Bw_{2m+1}-(m+1)bw_{2m+3}=0$ and
 fixes the normalisation:
$Ch_*(u)=(c_{2m+1}w_{2m+1})$ where 
$$w_{2m+1}=(u^{-1},u,\ldots,u^{-1},u)\ \ (2m+2\   \mbox{entries}),\ c_{2m+1}=
(-1)^{m}m!.$$
Now, in the normalised complex:
$$b(1,u^{-1},u,\ldots,u^{-1},u)=(u^{-1},u,\ldots,u^{-1},u)
+(u,u^{-1}\ldots,u,u^{-1})$$
and
$$B(1,u^{-1},u,\ldots,u^{-1},u)=0.$$
Hence let $z=(z_{2m+2})$ where $m=0,1,2,\ldots$
and 
$$z_{2m+2}=c_{2m+1}(1,u^{-1},u,\ldots,u^{-1}, u) \ \ (2m+3\ \mbox{ entries)}.$$

Then
$$Ch_{2m+1}(u^*)+Ch_{2m+1}(u)=c_{2m+1}b(1,u^{-1},u,\ldots,u^{-1},u)
+c_{2m-1}B(1,u^{-1},u,\ldots,u^{-1},u))$$
so that $Ch_*(u)+Ch_*(u^{-1})=(b+B)z$ 
\end{proof}

\section{The Main Result and Outline of the Proof}\label{sec4}

\subsection{Statement of the Main Result}

We introduce some notation in order to be able to
state the main theorem.

First, we require multi-indices $(k_1,...,k_m)$, $k_i\in\{0,1,2,...\}$, whose
length $m$ will always be
clear from the context. We write $|k|=k_1+\cdots+k_m$, and define
$\alpha(k)$ by
\ben \alpha(k)=\frac{1}{k_1!k_2!\cdots
k_m!(k_1+1)(k_1+k_2+2)\cdots(|k|+m)}.\een
The numbers $\s_{n,j}$ are defined by the equality
\ben\prod_{j=0}^{n-1}(z+j+1/2)=\sum_{j=0}^{n}z^j\s_{n,j}.\een
These are just the elementary symmetric functions of $1/2,3/2,...,n-1/2$.

If $(\A,\HH,\D)$ is a $QC^\infty$ spectral triple and $T\in\cn$, we write
$T^{(n)}$ to denote  the iterated commutator
$[\D^2,[\D^2,[\cdots,[\D^2,T]\cdots]]]$ where we have $n$ commutators with
$\D^2$. It follows from the remarks after Definition \ref{qck} that
operators of the form $T_1^{(n_1)}\cdots
T_k^{(n_k)}(1+\D^2)^{-(n_1+\cdots+n_k)/2}$ are in $\cn$ when
$T_i=[\D,a_i]$, or $=a_i$ for $a_i\in\A$.

\begin{defn}\label{dimension}
If $(\A,\HH,\D)$ is a $QC^\infty$ spectral triple,
we call
\ben p=\inf\{k\in{\R}:\tau((1+\D^2)^{-k/2})<\infty\}\een
the {\bf spectral dimension} of $(\A,\HH,\D)$.
We say that $(\A,\HH,\D)$ has {\bf isolated spectral dimension} if
for $b$ of the form
$$b=a_0[\D,a_1]^{(k_1)}\cdots[\D,a_m]^{(k_m)}(1+\D^2)^{-m/2-|k|}$$
the zeta functions
\ben \zeta_b(z-(1-p)/2)=\tau(b(1+\D^2)^{-z+(1-p)/2})\een
have analytic continuations to a deleted neighbourhood of $z=(1-p)/2$.
\end{defn}

{\bf Remark} Observe that we allow the possibility that the analytic 
continuations of these zeta functions may have an essential singularity at 
$z=(1-p)/2$. All that is necessary for us is that the residues at this point 
exist.

Now we define, for $(\A,\HH,\D)$ having isolated spectral dimension and
$$b=a_0[\D,a_1]^{(k_1)}\cdots[\D,a_m]^{(k_m)}(1+\D^2)^{-m/2-|k|}$$
\ben \tau_j(b)=res_{z=(1-p)/2}(z-(1-p)/2)^j\zeta_b(z-(1-p)/2).\een
The hypothesis of isolated spectral dimension is clearly necessary here in 
order to define the residues.

With these preliminaries we can state the main result of the paper.

\begin{thm}[Semifinite Odd Local Index Theorem]\label{SFLIT} 
Let $(\A,\HH,\D)$ be an 
odd finitely summable $QC^\infty$ spectral triple with spectral
dimension $p\geq 1$. Let $N=[p/2]+1$ where $[\cdot]$ denotes the 
integer part, and let $u\in\A$ be unitary. Then

1) \qquad\qquad\qquad $sf(\D,u^*\D u)=\frac{1}{\sqrt{2\pi i}} 
res_{r=(1-p)/2}\left( 
\sum_{m=1,odd}^{2N-1} \phi_m^r(Ch_m(u))\right)$

where for $a_0,...,a_m\in\A$, $l=\{a+iv:v\in{\R}\}$, $0<a<1/2$, 
$R_s(\lambda)=(\lambda-(1+s^2+\D^2))^{-1}$ and $r>0$ we define 
$\phi_m^r(a_0,a_1,...,a_m)$ to be
\ben \frac{-2\sqrt{2\pi i}}{\Gamma((m+1)/2)}
\int_0^\infty s^m\tau\left(\frac{1}{2\pi i}
\int_l\lambda^{-p/2-r}a_0R_s(\lambda)[\D,a_1]R_s(\lambda)\cdots
[\D,a_m]R_s(\lambda)d\lambda\right)ds.\een
In particular the sum on the right hand side of $1)$ analytically continues 
to a deleted neighbourhood of $r=(1-p)/2$ with {\em at worst} a simple pole 
at $r=(1-p)/2$.
Moreover, the complex function-valued cochain $(\phi_m^r)_{m=1,odd}^{2N-1}$ 
is a $(b,B)$ 
cocycle for $\A$ modulo functions holomorphic in a half-plane containing 
$r=(1-p)/2$.

2) The spectral flow  $sf(\D,u^*\D u)$ is also the residue of a sum of zeta 
functions:
\bean &&\frac{1}{\sqrt{2\pi i}} res_{r=(1-p)/2} \left(\sum_{m=1,odd}^{2N-1}
\sum_{|k|=0}^{2N-1-m}\sum_{j=0}^{|k|+(m-1)/2}(-1)^{|k|+m}
\alpha(k)\Gamma((m+1)/2)\s_{|k|+(m-1)/2,j} \right.\nno
&& \qquad\qquad\biggl.(r-(1-p)/2)^j\tau\left(u^*[\D,u]^{(k_1)}
[\D,u^*]^{(k_2)}\cdots[\D,u]^{(k_m)}(1+\D^2)^{-m/2-|k|-r+(1-p)/2}\right)
\Biggr).\eean
In particular the sum of zeta functions on the right hand side analytically 
continues to a deleted neighbourhood of $r=(1-p)/2$ and has {\rm at worst} a 
simple pole at 
$r=(1-p)/2$.

3) If $(\A,\HH,\D)$ also has isolated spectral dimension then
\ben sf(\D,u^*\D u)=\frac{1}{\sqrt{2\pi i}}\sum_m \phi_m(Ch_m(u))\een
where for $a_0,...,a_m\in\A$
\bean \phi_m(a_0,...,a_m)&=&res_{r=(1-p)/2}\phi^r_m(a_0,...,a_m)=\sqrt{2\pi i}
\sum_{|k|=0}^{2N-1-m}(-1)^{|k|}
\alpha(k)\times\nno
&\times&\sum_{j=0}^{|k|+(m-1)/2}\s_{(|k|+(m-1)/2),j}\tau_j
\left(a_0[\D,a_1]^{(k_
1)}\cdots[\D,a_m]^{(k_m)}(1+\D^2)^{-|k|-m/2}\right),\eean
and $(\phi_m)_{m=1,odd}^{2N-1}$ is a $(b,B)$ cocycle for $\A$. When $[p]=2n$ 
is even, the term with $m=2N-1$ is zero, and 
for $m=1,3,...,2N-3$, all the top terms with $|k|=2N-1-m$ are zero. 
\end{thm}

{\bf Remark} Since $\phi_m$ is a multilinear functional, it is well-defined on 
elements of $\A^{\otimes m+1}$ such as $Ch_m(u)$.

\begin{cor}\label{lowdim}
For $1\leq p<2$, the statements in $3)$ of Theorem \ref{SFLIT} are true without 
the 
assumption of isolated dimension spectrum.
\end{cor}

\subsection{Outline of the Proof of Theorem \ref{SFLIT}}

The basic strategy we employ is very simple, but as the technical details may 
obscure this, we offer an outline of the proof here using the notation of 
Theorem \ref{SFLIT}. There are two basic parts of the proof.

{\bf 1.} 
The Carey-Phillips spectral flow formula is manipulated into a form that allows
us to use perturbation theory in the form of a
resolvent expansion. The resulting formula suggests the definition of a 
substitute, for finitely summable spectral triples, of the JLO cocycle of 
entire cyclic cohomology. 
Our substitute we term the `resolvent cocycle'. It
is a function-valued $(b,B)$-cocycle, 
modulo functions holomorphic in a certain half-plane. 

{\bf 2.} The pseudodifferential calculus of \cite{CM} then enables us 
to write the spectral flow as a sum of zeta functions, 
modulo functions holomorphic in a certain half-plane.
If we impose the isolated spectral dimension assumption we can 
analytically continue these zeta functions and 
take residues at a predetermined critical point. 
We then see that spectral flow is obtained by pairing $Ch_*(u)$ with 
a variant of the Connes-Moscovici residue cocycle.

We now expand on these two basic parts.

In Sections \ref{nandt} and \ref{tande}, we prove all the basic norm and 
trace estimates we will require for the resolvent expansion. 
These technical details may be skipped on a first reading.

The proof proper begins in Section \ref{G}. To successfully apply 
perturbation techniques to the Carey-Phillips spectral flow formula, 
Equation (\ref{basicformula}), we require `more room to manoeuvre'. 
Three basic steps are involved in this section. First, we `double-up' the 
data $(\HH,\D)$ from our spectral triple and unitary $u$ 
by tensoring on two copies of ${\C}^2$ to $\mathcal H$ to obtain 
an unbounded self-adjoint
operator $\tilde\D$ affiliated with $\tilde{\mathcal N}:=
M_2\otimes M_2\otimes\mathcal N$ and a self-adjoint unitary $q$ 
(determined by $u$) in $M_2\otimes M_2\otimes\mathcal A$. 
This may be viewed as employing a formal (Clifford) Bott periodicity 
and replaces the trace $\tau$ by a supertrace $S\tau$. 

Second, with the additional freedom, we
 are now able use an idea of \cite{G}
to define a two parameter family of perturbations
$\tilde\D_{r,s}$, 
$r\in[0,1]$ and $s\in[0,\infty)$. We observe that, 
as the spectral flow is computed by integrating 
an {\em exact} one-form on an affine space of perturbations of $\tilde\D$, 
Lemma \ref{exact}, we may compute 
spectral flow from $\D$ to $u^*\D u$ along different paths 
joining the endpoints; initially it is given 
by integrating with respect to $r$ when $s=0$. 

The third step chooses a path that expresses spectral flow in terms of
an integral over the $s$ variable with $r=0$
where the perturbation
in (\ref{basicformula}), instead of being of first order in $\tilde\D$, 
is now zeroth order.
Thus we obtain a new formula for spectral flow
\ben sf(\D,u^*\D u)=\frac{1}{C_{p/2+r}}\int_0^\infty 
S\tau\left(q(1+\tD^2+s\{\tD,q\}+s^2)^{-p/2-r}\right)ds\een
where $\{\cdot,\cdot\}$ denotes the anticommutator. 
Crucially, the anticommutator $\{\tilde\D,q\}$ is bounded, and 
we are now in a position to employ perturbation theory in the form of the  
resolvent expansion.

Section \ref{psido} reviews the pseudodifferential calculus of \cite{CM}, 
and the 
`Taylor expansion' in the form introduced by Higson, \cite{hig}. We then prove 
several technical results, Lemmas \ref{normtrick}-\ref{uni}, that allow us to 
easily apply the pseudodifferential calculus in our setting. Again, this 
section may be omitted on a first reading.

In Section \ref{resexpand} we write
\ben (1+\tD^2+s^2+s\{\tD,q\})^{-p/2-r}=
\frac{1}{2\pi i}\int_l\lambda^{-p/2-r}
(\lambda-(1+\tD^2+s^2+s\{\tD,q\}))^{-1}d\lambda,\een
where the vertical line $l$ lies between $0$ and $spec(1+\tD^2+s^2+s\{\tD,q\})$
for all $s\in [0,\infty)$. We then apply the resolvent expansion 
(writing $R_s(\lambda)=(\lambda-(1+\tD^2+s^2))^{-1}$)
\ben (\lambda-(1+\tD^2+s^2+s\{\tD,q\}))^{-1}=
\sum_{m=0}^{2N-1}\left(R_s(\lambda)s\{\tD,q\}\right)^mR_s(\lambda)
+Remainder.\een
The estimates in Lemmas \ref{yucko} and \ref{spaz} allow us to show in 
Lemma \ref{resolvent} that {\em modulo functions of $r$ holomorphic in a 
half-plane containing} $r=(1-p)/2$
\be sf(\D,u^*\D u)C_{p/2+r}=
\frac{1}{2\pi i}\sum_{m=1,odd}^{2N-1}
\int_0^\infty s^mS\tau\left(q\int_l\lambda^{-p/2-r}
(R_s(\lambda)\{\tD,q\})^mR_s(\lambda)d\lambda\right) ds.\label{resexpeqn}\ee
The even terms in the above expansion are seen to 
vanish by elementary Clifford-type manipulations.
The `constant' 
\ben C_{p/2+r}=\frac{\Gamma(r-(1-p)/2)\Gamma(1/2)}{\Gamma(p/2+r)}\een
has simple poles at $r=(1-p)/2-k$, $k=0,1,2,...$, with residue equal to 
$1$ at $r=(1-p)/2$. Therefore, since the error terms in 
Equation (\ref{resexpeqn}) 
are holomorphic at $r=(1-p)/2$, we may take residues at $r=(1-p)/2$
of the analytic continuations of both sides of (\ref{resexpeqn})
even though the individual terms in this expansion need not analytically
continue.

Section \ref{rescocycle} begins by performing the `super' part of the trace 
to obtain a formula for the spectral flow in terms of the original spectral 
triple $(\A,\HH,\D)$ and the unitary $u$. The general structure of this 
formula suggests the definition of a function-valued $(b,B)$-cochain on the 
algebra $\A$. We call this the resolvent cocycle, and using techniques 
inspired by Higson, \cite{hig}, we show that this is a cocycle 
{\em modulo functions of $r$ holomorphic in a half-plane containing 
$(1-p)/2$}. This `almost cocycle' property is proved in Lemma \ref{cocycle}, 
and this proves $1)$ of Theorem \ref{SFLIT}.

Section \ref{residue} returns to our spectral flow computations. Section 
\ref{psiexpan} applies the pseudodifferential calculus, in the form derived 
in Lemma \ref{firstexpan}, to each term of the resolvent expansion. 
This moves all the resolvents to the right, allowing us to use Cauchy's 
formula to perform the complex line integral. In Section \ref{sint} we 
perform the remaining integral over $s\in[0,\infty)$, and so obtain our 
penultimate formula:
\bea &&sf(\D,u^*\D u)C_{p/2+r}\nno
&&\qquad=\sum_{m=1,odd}^{2N-1}\sum_{|k|=0}^{2N-1-m}
C_{k,m,r}S\tau\left(q\{\tD,q\}^{(k_1)}\cdots
\{\tD,q\}^{(k_m)}(1+\tD^2)^{-(p-1)/2-|k|-m/2-r}\right),
\label{remarkable}\eea
where equality is again {\em modulo functions of $r$ holomorphic in a 
half-plane containing $(1-p)/2$}.  

This is a remarkable formula. So far we have not invoked the isolated spectral 
dimension hypothesis, yet the sum of zeta functions in 
Equation \ref{remarkable} clearly has a simple pole at $r=(1-p)/2$, with 
residue equal to the spectral flow. This proves part $2)$ of 
Theorem \ref{SFLIT}.

In Section \ref{residuecocycle} we finally assume that the individual zeta 
functions possess analytic continuations to a deleted neighbourhod of 
$r=(1-p)/2$ so we can take residues of the zeta functions in Theorem 
\ref{grandfinale} to obtain our version of  the residue cocycle.
We can then prove part $3)$ of  Theorem \ref{SFLIT}. 
The cocycle property for the residue cocycle follows from the `almost' 
cocycle property of the resolvent cocycle upon taking residues. We conclude 
with the simple proof of Corollary \ref{lowdim}.

\section{Norm and Trace Estimates}
Throughout this section, let $D:\mbox{dom}D\subseteq\HH\tto\HH$ 
be an unbounded self-adjoint operator on the Hilbert space $\HH$.
\subsection{Norm Estimates}\label{nandt}

In a number of estimates, we will also consider a bounded self-adjoint operator
$A$. The operators $A$ that are of interest satisfy $s^2+sA+D^2\geq 0$ for 
all real $s\geq 0$. However it is also convenient at times to assume that
$\n A\n$ is relatively small: $\n A\n<\sqrt{2}$, for example. This can be 
achieved by scaling $A$: see {\bf Observation 2} of section 7.

\begin{lemma}\label{normlambda} Let $D$ be an unbounded self-adjoint 
operator.\\
(a) For $\lambda=a+iv\in \C$, $0<a<1/2$, $s\geq 0$ we have the estimate 
\ben \n(\lambda-(1+D^2+s^2))^{-1}\n\leq (v^2+(1+s^2-a)^2)^{-1/2}
\leq\frac{1}{1-a}.\een
(b) If $A$ is bounded, self-adjoint and $s^2+sA+D^2\geq 0$ we have 
\ben \n(\lambda-(1+D^2+s^2+sA))^{-1}\n\leq 
(v^2+(1-a)^2)^{-1/2}\leq\frac{1}{1-a}.\een
(c) If $A$ is bounded, self-adjoint and $c=\n A\n<\sqrt{2}$ we have
\ben \n(\lambda-(1+D^2+s^2+sA))^{-1}\n\leq 
(v^2+(1+s^2-a-sc)^2)^{-1/2}\leq\frac{1}{1/2-a}.\een
\end{lemma}
\begin{proof} Part (a) is an application of the functional 
calculus. To see (b) note that the spectrum of 
$(1+D^2+s^2+sA)$ is contained in $[1,\infty)$, so the distance from 
$\lambda$ to the spectrum is at least $|\lambda -1|$. For (c)
note that the minimum value in the spectrum of
$(1+D^2+s^2+sA)$ is at least $(1+s^2-sc)>1/2$, and so the distance from
$\lambda$ to the spectrum is at least $|\lambda -(1+s^2-sc)|.$
\end{proof}
 
\subsection{Trace Estimates}\label{tande}

For the duration of this section we suppose that the operator $D$ satisfies 
the 
summability condition
\be (1+D^2)^{-n/2}\in \LL^1({\mathcal N})\ \ \ \forall n>p\geq 1.
\label{summ}\ee
For instance, if $\D$ comes from a spectral triple $(\A,\HH,\D)$ then this is
exactly the condition that $\D$ have dimension $p$. 
With the hypothesis (\ref{summ}) we will obtain some trace norm estimates for 
various operators depending on $s$. The following Lemma is the key technical 
estimate, and the main  result of this subsection.
Henceforth $1>>\epsilon>0$. 

\begin{rems*}
In the following lemma, if instead of assuming  $\n A\n<\sqrt{2}$, we suppose
that\\ $D^2+s^2+sA\geq 0$, then the estimate on the RHS becomes
$C_{p+\epsilon}(1/2+(1/2)s^2)^{-Re(r)+\epsilon}$ for $s\geq 2\n A\n$. Since 
one can show that $(1+D^2+s^2+sA)^{-p/2-r}$ is trace-class for all s, the
integrability conclusion still holds also.
\end{rems*}

\begin{lemma}\label{BIG} Let $A$ be bounded self-adjoint operator
with $\n A\n<\sqrt{2}$. Let $p\geq 1$ and let $(1+D^2)^{-1/2}$
be $(p+\epsilon)$-summable for every $\epsilon >0$. Then for each $\epsilon >0$
and $r\in{\C}$ with $Re(r)>0$, the trace norm of $(1+D^2+s^2+sA)^{-p/2-r}$ 
satisfies
$$\n (1+D^2+s^2+sA)^{-p/2-r}\n_1 \leq 
C_{p+\epsilon}(1/2+s^2-s\n A\n)^{-Re(r)+\epsilon},$$
where $C_{p+\epsilon}= \n (1/2+D^2)^{-(p/2+\epsilon)}\n_1.$
So, if $Re(r)>1/2+\epsilon,$ then $\n (1+D^2+s^2+sA)^{-p/2-r}\n_1$ 
is integrable in s on ${\bf R}$.
\end{lemma}

\begin{proof}
Throughout the proof we suppose, without loss of generality, 
that $r$ is real and positive. 
This is possible because for complex $r$ with $Re(r)>0$ we have
\bean &&\n (1+D^2+s^2+sA)^{-p/2-r}\n_1\nno
&\leq& 
\n(1+D^2+s^2+sA)^{-iIm(r)}\n\n(1+D^2+s^2+sA)^{-p/2-Re(r)}\n_1\nno
&\leq &\n(1+D^2+s^2+sA)^{-p/2-Re(r)}\n_1.\eean
 We first consider the case where $A=0$. 
Since $(1/2+D^2)^{-1}\leq 2(1+D^2)^{-1}$ the constant $C_{p+\epsilon}$
is finite.
Now, for positive real numbers $X,Y,a,b$ one easily sees that:
$$(X+Y)^{-a}\leq X^{-a}\ \ \mbox{and}\ \ (X+Y)^{-b}\leq Y^{-b}$$
and so,
$$(X+Y)^{-a-b}\leq X^{-a}Y^{-b}.$$ 
Therefore, by the functional calculus:
\begin{eqnarray*}
(1+D^2+s^2)^{-p/2-r}&=&(1/2+D^2+1/2+s^2)^{-(p/2+\epsilon)-(r-\epsilon)}\\
&\leq& (1/2+D^2)^{-(p/2+\epsilon)}(1/2+s^2)^{-(r-\epsilon)}.
\end{eqnarray*}
Taking the trace norm of this inequality, gives the lemma when $A=0$.

To obtain the general case where $\n A\n<\sqrt{2}$, we observe that
\ben 0<1+s^2+D^2-s\n A\n\leq 1+s^2+D^2+sA\leq 1+s^2+D^2+s\n A\n.\een
Consequently, 
\ben (1+s^2+D^2+sA)^{-1}\leq(1+s^2+D^2-s\n A\n)^{-1}.\een
Now Corollary 4, Appendix B of \cite{CP1} tells us that
\ben \tau((1+s^2+D^2+sA)^{-p/2-r})\leq
\tau((1+s^2+D^2-s\n A\n)^{-p/2-r}).\een
By the same argument as the case $A=0$ (noting that $1/2+s^2-s\n A\n >0$) 
we get:
$$(1+s^2+D^2-s\n A\n)^{-p/2-r}\leq (1/2+D^2)^{-(p/2+\epsilon)}
(1/2+s^2-s\n A\n)^{-(r-\epsilon)}).$$
Taking the trace of this inequality and combining it with the immediately 
previous inequality, yields the proof of the lemma.
\end{proof}

\begin{lemma}\label{lambda} Let $0<a=Re(\lambda)<1/2$, $\lambda=a+iv$. 
Let $p\geq 1$ and let $(1+D^2)^{-1/2}$
be $(p+\epsilon)$-summable for every $\epsilon >0$. 
Then for each $\epsilon >0$ and 
$N>(p+\epsilon)/2$ , we have the trace-norm estimate:  
$$\n (\lambda-(1+D^2+s^2))^{-N}\n_1\leq C'_{p+\epsilon}
((1/2+s^2-a)^2+v^2)^{-N/2+(p+\epsilon)/4}$$
where $C'_{p+\epsilon}= \n (1/2+D^2)^{-(p+\epsilon)/2}\n_1.$
\end{lemma}

\begin{proof} This is similar to Lemma \ref{BIG}. 
Now,
\begin{eqnarray*}
|\left[\lambda -(1+D^2 +s^2)\right]^{-1}|&=&
\left[(1/2 + D^2 + 1/2 +s^2-a)^2+v^2\right]^{-1/2}\\
&\leq&\left[(1/2 + D^2)^2 +(1/2+s^2-a)^2+v^2\right]^{-1/2}
\end{eqnarray*}
By the argument of the previous lemma, we have:
\begin{eqnarray*}
|\left[\lambda -(1+D^2 +s^2)\right]^{-N}|
&\leq&\left[(1/2 + D^2)^2 +(1/2+s^2-a)^2+v^2\right]^{-N/2}\\
&=&\left[(1/2 + D^2)^2 +(1/2+s^2-a)^2+v^2\right]
^{-(p+\epsilon)/4-(N/2-(p+\epsilon)/4)}\\
&\leq& (1/2 + D^2)^{-(p+\epsilon)/2}((1/2+s^2-a)^2+v^2)^{-N/2+(p+\epsilon)/4}.
\end{eqnarray*}
Taking the trace of this inequality yields the proof of the lemma.
\end{proof}

We finish this 
subsection with an integral estimate which we will use several times.

\begin{lemma}\label{intest} Let $0<a<1/2$ and $0\leq c\leq\sqrt{2}$ 
and $A=0$ or $1$. Let $J$,$K$, and $M$ be nonegative constants. 
Then the integral
\be \int_0^\infty \int_{-\infty}^\infty s^J\sqrt{a^2+v^2}^{-M}
\sqrt{(s^2+1/2-a)^2+v^2}^{-K}\sqrt{(s^2+1-a-sc)^2+v^2}^{-A}dvds
\label{integral}\ee
converges provided $J-2K-2A<-1$ and $J-2K-2A+1-2M<-2$.
\end{lemma}

\begin{proof} We begin with the case $A=0$. The integrand is positive and 
continuous (and hence measurable), so by 
Tonelli's theorem we may evaluate the $s$ integral first. As $1/2-a>0$,
\ben (s^2+1/2-a)^2\geq s^4+(1/2-a)^2\een
and
\ben \int_0^\infty s^J((s^2+1/2-a)^2+v^2)^{-K/2}ds\leq\int_0^\infty 
s^J(s^4+(1/2-a)^2+v^2)^{-K/2}ds.\een
Now write $b^2=(1/2-a)^2+v^2$ and set $t=sb^{-1/2}$. Performing the substitution
 gives
\bean \int_0^\infty 
s^J(s^4+b^2)^{-K/2}ds&=&b^{-K}\int_0^\infty s^J(s^4/b^2 +1)^{-K/2}ds\nno
&=& b^{-K+J/2+1/2}\int_0^\infty t^J(t^4+1)^{-K/2}dt.\eean
This integral converges provided $J-2K<-1$. Thus 
\bean &&\int_0^\infty s^J\int_{-\infty}^\infty\sqrt{a^2+v^2}^{-M}
\sqrt{(s^2+1/2-a)^2+v^2}^{-K}dvds\nno
&&\leq C\int_{-\infty}^\infty 
((1/2-a)^2+v^2)^{-K/2+J/4+1/4}(a^2+v^2)^{-M/2}dv,\eean
and this is finite provided $-2K+J+1-2M<-2$. When $A=1$ we observe that
\ben 1+s^2-a-sc=(s-c/2)^2+(1-c^2/4-a)\een
so
\bean (1+s^2-a-sc)^2&=&((s-c/2)^2+(1-c^2/4-a))^2\nno
&\geq & (s-c/2)^4+(1-c^2/4-a)^2\nno
&\geq & (s-c/2)^4+(1/2-a)^2.\eean
We also have for $s\geq c/2$
\ben (s^2+1/2-a)^2\geq (s-c/2)^4+(1/2-a)^2.\een
For $s\leq c/2$ we have 
\ben (s^2+1/2-a)^2\geq (1/2-a)^2.\een
The integral in Equation (\ref{integral}) is then bounded by
\bean &&\int_{c/2}^\infty s^J\int_{-\infty}^\infty\sqrt{a^2+v^2}^{-M}
\sqrt{(s-c/2)^4+(1/2-a)^2+v^2}^{-K-1}dvds\nno
&+& \int_0^{c/2}(c/2)^J\int_{-\infty}^{\infty}\sqrt{a^2+v^2}^{-M}
\sqrt{(1/2-a)^2+v^2}^{-K-1}dvds.\eean
We look at the integral over $s$ from $c/2$ to $\infty$. Write 
$b^2=(1/2-a)^2+v^2$, so we are considering
\ben \int_{c/2}^\infty s^J\sqrt{(s-c/2)^4+b^2}^{-K-1}ds=
\int_0^\infty (s+c/2)^J\sqrt{s^4+b^2}^{-K-1}ds.\een
Now setting $t=sb^{-1/2}$ and making the change of variables this last integral 
becomes
$$\int_0^\infty 
b^{-K+(J-1)/2}(t+cb^{-1/2}/2)^J(1+t^4)^{-(K+1)/2}$$
\ben \leq b^{-K-1/2+J/2}\int_0^\infty (t+c/2\sqrt{1-2a})^J(1+t^4)^{-(K+1)/2}dt,
\een
and this converges when $J-2K-2<-1$. Now we have the $v$-integral, 
which is bounded by
\ben C\int_{-\infty}^\infty 
\sqrt{a^2+v^2}^{-M}\sqrt{(1/2-a)^2+v^2}^{-K-1/2+J/2}dv\een
and this is finite when $J-2K-1-2M<-2$. The  $s$-integral from 
$0$ to $c/2$ is finite, and the corresponding $v$ integral is given by
\ben \int_{-\infty}^\infty \sqrt{a^2+v^2}^{-M}\sqrt{(1/2-a)^2+v^2}^{-K-1}dv\een
and this is finite for $-M-K-1<-1$ or $M+K>0$.
\end{proof}

\subsection{Application: rewriting the formula for spectral flow.}\label{G}
In this subsection 
we begin with the spectral flow formula (\ref{basicformula}) of 
the last Section, for a finitely summable odd spectral triple $(\A,\HH,\D)$
and rewrite it in a different way so as to be able to exploit resolvent
expansions. 

Our method borrows an idea from \cite{G} however, whereas Getzler's
approach
is via the superconnection formalism, we will adopt a more concrete
functional analytic approach suggested by
\cite{CP0}. See also section 9 of \cite{CP2}.

\begin{defn}
Form the Hilbert space
$\tilde{\mathcal H}=\IC^2\otimes\IC^2\otimes\mathcal H$ acted on by the
von Neumann algebra,
$\tilde{\mathcal N}=M_2\otimes M_2\otimes\mathcal N$.
Introduce the two dimensional Clifford algebra in the form
$$\sigma_1= \left(\begin{array}{cc}
                   0 & 1 \\
                   1 & 0
                   \end{array} \right),\ \  \sigma_2= \left(\begin{array}{cc}
                   0 & -i \\
                   i & 0
                   \end{array} \right), \ \  \sigma_3=
 \left(\begin{array}{cc} 1 & 0 \\
                         0 & -1
\end{array} \right).$$
(With
$\tD$ defined below this gives an odd spectral triple for the algebra 
$1_2\otimes M_2\otimes\mathcal A$ but we shall {\bf never} use this fact). 
Define the grading on $\tilde{\mathcal H}$ by
$\Gamma = \sigma_2\otimes\sigma_3\otimes 1\in \tilde{\mathcal N}$. 
\end{defn}
Let $u\in \mathcal A$ be unitary and
introduce the following even operators (i.e., they commute with $\Gamma$):
$$\tD = \sigma_2\otimes 1_2\otimes \D, \ \
q=\sigma_3\otimes\left(\begin{array}{cc}
                   0 & -iu^{-1} \\
                   iu & 0
                   \end{array} \right),\ \
\D_r= (1-r)\tD - rq\tD q, \ \   
\D_{r,s} = \D_r + sq$$
for $r\in [0,1], s\in [0,\infty).$
Clearly, the unbounded operators are affiliated
with $\tilde{\mathcal N}$. Notice that
$$ \D_{r}\equiv \D_{r,0}= \sigma_2\otimes \left(\begin{array}{cc}
                   \D+ru^{-1}[\D,u] & 0 \\
                   0 & \D+ru[\D,u^{-1}]
                   \end{array} \right).$$
So 
$$\dot \D_r = \sigma_2\otimes \left(\begin{array}{cc}
                   u^{-1}[\D,u] & 0 \\
                   0 & u[\D,u^{-1}]
                   \end{array} \right).$$
The graded trace on $\tilde{\mathcal N}$ we write as
$S\tau(a) =\frac{1}{2}\tau(\Gamma a)$ for $a$ trace-class in 
$\tilde{\mathcal N}$, and so for example,
\bean &&\qquad\qquad\qquad\qquad\qquad 
S\tau(\dot \D_r(1+\D_r^2)^{-\frac{n}{2}})\nno
&=&\frac{1}{2}\int_0^1\tau\left(1_2\otimes\bma 
u^*[\D,u](1+\D_r^2)^{-n/2} & 0\\ 0 & 
-u[\D,u^*](1+\D_r^2)^{-n/2}\ema\right)dr\nno
 &=&\tau\left(u^{-1}[\D,u](1+(\D+ru^{-1}[\D,u])^2)^{-\frac{n}{2}}- u[\D,u^{-1}] 
(1+(\D+ru[\D,u^{-1}])^2)^{-\frac{n}{2}}\right).\eean
Next we calculate
$$\D_{r,s}^2= \D_r^2+ s(1-2r)\sigma_1\otimes\left(\begin{array}{cc}
                   0 & [\D,u^{-1}] \\
                   -[\D,u] & 0
                   \end{array} \right) +s^2$$
which depends on the relation
$$\tD q + q\tD =\sigma_1\otimes\left(\begin{array}{cc}
                   0 & [\D,u^{-1}] \\
                   -[\D,u] & 0
                   \end{array} \right) $$
Our immediate aim now is to rewrite the spectral flow in terms of these new 
operators. 
\begin{lemma} \label{exact}
Let $(\A,\HH,\D)$ be a spectral triple with dimension $p\geq 1$. Consider
the affine space $\Phi$ of perturbations of $\tilde\D$
given by
$$\Phi=\{\tilde\D+X \ \ |\  X\in\tilde{\mathcal N}_{sa}\; \mbox{and}\
[X,\Gamma]=0\}.$$
Then for $n>p$
the linear functional on the tangent space to $\Phi$ at $\tD_1\in\Phi$
$$\alpha_{\tD_1}(X)=S\tau(X(1+\tD_1^2)^{-n/2})$$
is well defined and makes
$\tD_1\mapsto \alpha_{\tD_1}$ an exact one-form (ie an exact section of the
cotangent bundle to $\Phi$).
\end{lemma}

\begin{proof}
We would like to apply Theorem 9.3 of \cite{CP2} to the pair 
$(\tilde\cn,\tD)$; however, that theorem 
deals with the ungraded case (no $\Gamma$ and just the ordinary trace $\tau$) 
and this lemma is
firmly in the graded setup. In order to avoid reworking the whole theory
to fit the graded setting, we resort to the following trick. We modify our
unbounded operator, setting $\hat{\D}:=\Gamma\tD$, and apply Theorem 9.3
of \cite{CP2} to the pair $(\tilde\cn,\hat{D})$. Since $\Gamma$ is 
self-adjoint and commutes with all the
operators in $\Phi$, we see that $\Gamma\Phi$ is a real affine space and
$$\Gamma\Phi=\{\hat{\D}+X \ \ |\  X\in\tilde{\mathcal N}_{sa}\; \mbox{and}\
[X,\Gamma]=0\}\subseteq \hat{\D} + \tilde{\mathcal N}_{sa}.$$
Thus, if $\tD+X(s)$ is a piecewise $C^1$ continuous path in $\Phi$, 
then $\hat{\D}+\Gamma X(s)$ is a piecewise $C^1$ continuous path in 
$\Gamma\Phi\subseteq \hat{\D} + \tilde{\mathcal N}_{sa}$,
and we have by Theorem 9.3 of \cite{CP2} that the integral:
$$\int_0^1 \tau\left(\frac{d}{ds}\left(\Gamma X(s)\right)
\left(1+(\hat{\D}+\Gamma X(s))^2\right)^{-n/2}\right)ds$$
depends only on the endpoints $\Gamma X(0)$ and $\Gamma X(1)$ and hence only
on $X(0)$ and $X(1)$. But, since
$$(\hat{\D} +\Gamma X(s))^2=(\Gamma[\tD+X(s)])^2=(\tD+X(s))^2$$
we see that our integral becomes:
$$\int_0^1\tau\left(\Gamma X^{\prime}(s)\left(1+(\tD+X(s))^2\right)^{-n/2}
\right)ds=2\int_0^1 S\tau\left(X^{\prime}(s)
\left(1+(\tD+X(s))^2\right)^{-n/2}\right)ds,$$
and its value depends only on the endpoints $X(0)$ and $X(1)$. Since integrals
of our one-form are independent of path, our one-form is exact by Lemma 7.3
of \cite{CP2}. In fact, it is this independence of path that we use below.
\end{proof}

\begin{lemma}\label{Lemma4.1} For $(\A,\HH,\D)$ a spectral triple with dimension
$p\geq 1$ we have for $n>p$
\ben \frac{1}{C_{n/2}}\int_0^1S\tau(\dot{\D_r}(1+\D_r^2)^{-n/2})dr=
2sf(\D,u^*\D u).\een
\end{lemma}

\begin{proof} Using the above computations we have 
\bean &&\int_0^1S\tau(\dot{\D_r}(1+\D_r^2)^{-n/2})dr\nno
&=&\int_0^1\tau\left( u^*[\D,u](1+(\D+ru^*[\D,u])^2)^{-n/2} 
-u[\D,u^*] (1+(\D+ru[\D,u^*])^2)^{-n/2}\right)dr\nno
&=&2{C_{n/2}} sf(\D,u^*\D u).\eean
By formula $(2)$ of section 4.1.
\end{proof}

\noindent{\bf Remark} Observe that this is the spectral flow from $\D$ to 
$u^*\D u$, which 
is $-Index(PuP)$ where $P$ is the spectral projection of $\D$ corresponding to 
the 
half-line $[0,\infty)$. This in turn is $-sf(\D,u\D u^*)$.
The following estimate enables us to exploit exactness of our one form
by changing the path of integration.

\begin{lemma}\label{Lemma3} Let $(\A,\HH,\D)$ be a spectral triple of dimension 
$p\geq 1$.
Then if $n>p$, we have
\ben 
\lim_{s_0\to\infty}\int_0^1S\tau\left(\frac{d\D_{r,s_0}}{dr}
(1+\D_{r,s_0}^2)^{-n/2
}\right)dr=0.\een
\end{lemma}

\begin{proof} 
First observe that from the definitions we have
\ben 
S\tau(\frac{d\D_{r,s_0}}{dr}(1+\D_{r,s_0}^2)^{-n/2})=
S\tau(\frac{d\D_{r,0}}{dr}(1+\D_{r,s_0}^2)^{-n/2}).\een
Next $\D_{1/2}=(1/2)(\tD-q\tD q)$, so that $\D_{1/2}$ anticommutes with $q$. 
Hence 
$(\D_{1/2}+s_0q)^2=\D_{1/2}^2+s_0^2$. So
\ben 
S\tau(\dot \D_{1/2}(1+\D_{1/2,s_0}^2)^{-n/2})
=S\tau(\dot \D_{1/2}(1+\D_{1/2}^2+s_0^2)^{-n/2}).\een
By Lemma \ref{BIG} with $r=(n-p)/2$, we have the following estimate 
(for a positive constant $C_{p+\epsilon}$)
$$\n(1+\D_{1/2}^2+s_0^2)^{-n/2}\n_1\leq 
C_{p+\epsilon}(1/2+s_0^2)^{-n/2+p/2+\epsilon}.$$
Since this exponent is negative (if we choose $\epsilon$ sufficiently small so 
that $n>p+2\epsilon$), we see that
as $s_0\to\infty$,
\bean S\tau(\dot\D_{1/2}(1+\D_{1/2,s_0}^2)^{-n/2})&\leq &S\tau 
(\left|\dot\D_{1/2}(1+\D_{1/2,s_0}^2)^{-n/2}\right|)\nno
&\leq & C_{p+\epsilon}(1/2+s_0^2)^{-n/2+p/2+\epsilon}\to 0.\eean
Set $A_r=\D_r-\D_{1/2}=(1/2-r)(\tD+q\tD q)$, and observe that there is a 
positive 
constant $c$ such that $\n A_r\n\leq c$ for all $r\in[0,1]$. Then 
\ben \D_{r,s_0}=\D_{1/2,s_0}+\D_{r,s_0}-\D_{1/2,s_0}=\D_{1/2,s_0}+A_r.\een
Using \cite[Corollary 8, page 710]{CP1} we obtain
\ben \n(1+\D_{r,s_0}^2)^{-n/2}\n_1\leq 
f(c)^{n/2}\n(1+\D_{1/2,s_0}^2)^{-n/2}\n_1,\een
where $f(c)=1+c^2/2+(c/2)\sqrt{c^2+4}$. Hence 
\ben \int_0^1 S\tau(\frac{d\D_{r,s_0}}{dr}(1+\D_{r,s_0}^2)^{-n/2})dr
\to 0\ \ \ \mbox{as}\ \ \ 
s_0\to\infty.\een
\end{proof}

\begin{lemma}\label{horizint} Let $(\A,\HH,\D)$ 
be a spectral triple of dimension 
$p\geq 1$. For $n>p$ and $s_0 > 0$ we have
\ben \int_0^{s_0} 
S\tau\left(\frac{d\D_{1,s}}{ds}(1+\D_{1,s}^2)^{-n/2}\right)ds=-\int_0^{s_0} 
S\tau\left(\frac{d\D_{0,s}}{ds}(1+\D_{0,s}^2)^{-n/2}\right)ds\een
\end{lemma}

\begin{proof}First observe that $\D_{1,s}=-q\tD q+sq$ and $\D_{0,s}=\tD +sq$
so that $$\frac{d\D_{1,s}}{ds}=q=\frac{d\D_{0,s}}{ds}.$$
Set $\rho=\s_2\otimes 1_2\otimes 1$, so that 
$\rho q\rho=-q$, $\rho^2=1$ and $\rho\Gamma\rho=\Gamma$. Then one easily 
calculates that:
$$\rho q\D_{1,s}=-(\tD +sq)\rho q = -\D_{0,s}\rho q.$$
So that:
$$\rho q\D_{1,s}^2=\D_{0,s}^2\rho q,$$
and hence for any Borel function, $f$:
$$\rho q f(\D_{1,s}^2)=f(\D_{0,s}^2)\rho q.$$
 Then
\bean 
S\tau(\frac{d\D_{1,s}}{ds}(1+\D_{1,s}^2)^{-n/2})&=&\frac{1}{2}\tau(\rho^2\Gamma
q(1+\D_{1,s}^2)^{-n/2})\nno
&=&\frac{1}{2}\tau(\Gamma\rho q(1+\D_{1,s}^2)^{-n/2}\rho)\nno
&=&\frac{1}{2}\tau(\Gamma (1+\D_{0,s}^2)^{-n/2}\rho q\rho)\nno
&=&-\frac{1}{2}\tau(\Gamma (1+\D_{0,s}^2)^{-n/2}q)\nno
&=&-\frac{1}{2}\tau(\Gamma q(1+\D_{0,s}^2)^{-n/2})\nno
&=&-S\tau(\frac{d\D_{0,s}}{ds}(1+\D_{0,s}^2)^{-n/2}).\eean
This completes the proof.
\end{proof}

\begin{lemma}\label{Corollary4.7} If $(\A,\HH,\D)$ is a spectral triple of 
dimension $p\geq 1$, then 
for $n>p$ we have
\ben sf(\D,u^*\D u)=\frac{1}{C_{n/2}}\int_0^\infty 
S\tau(\frac{d\D_{0,s}}{ds}(1+\D_{0,s}^2)^{-n/2})ds.\een
\end{lemma}

\begin{proof} The exactness of our one form gives us
\bean &&\frac{1}{C_{n/2}}\int_0^{s_0}
S\tau\left(\frac{d\D_{1,s}}{ds}(1+\D_{1,s}^2)^{-n/2}\right)ds-\frac{1}{C_{n/2}}
\int_0^{s_0} S\tau\left(\frac{d\D_{0,s}}{ds}(1+\D_{0,s}^2)^{-n/2}\right)ds\nno
&&+\frac{1}{C_{n/2}}\int_0^1S\tau\left(\frac{d\D_{r,0}}{dr}(1+\D_{r,0}^2)^{-n/2}
\right)dr=\frac{1}{C_{n/2}}\int_0^1S\tau\left(\frac{d\D_{r,s_0}}{dr}
(1+\D_{r,s_0}^2)^{-n/2}\right)dr\eean
Rearranging, using Lemma \ref{horizint} to combine the first two integrals and 
then
Lemma \ref{Lemma4.1} to substitute $2sf(\D,u^*\D u)$ for the third integral 
gives:
\ben sf(\D,u^*\D u)=\frac{1}{C_{n/2}}\int_0^{s_0}\hspace{-.1in}
S\tau\left(\frac{d\D_{0,s}}{ds}(1+\D_{0,s}^2)^{-n/2}\right)ds
+\frac{1}{2C_{n/2}}\int_0^1\hspace{-.1in}S\tau\left(\frac{d\D_{r,s_0}}{dr}
(1+\D_{r,s_0}^2)^{-n/2}\right)dr.\een
We can now take the limit as $s_0\to\infty$ of the right hand side using 
Lemma \ref{Lemma3} to
show that the improper integral converges to the spectral flow as claimed.
\end{proof}

The conclusion we draw from this is 
a shift of the path of integration of the integral
calculating spectral flow.
This leads us to a new formula
which, on setting, in Lemma \ref{Corollary4.7}, $n=p+2r$ we may write as:
\be sf(\D,u^*\D u)=\frac{1}{C_{p/2+r}}\int_0^\infty 
S\tau\left(q(1+\tD^2+s\{\tD,q\}+s^2)^{-p/2-r}\right)ds\label{uberformula}\ee
where $\{\cdot,\cdot\}$ denotes the anticommutator.
This formula will be the starting point for writing the spectral 
flow in terms of $(b,B)$ cocycles. 

Before discussing the pseudodifferential calculus, we wish to point out that 
while the improper integral in Equation (\ref{uberformula}) converges for 
$r>0$, it converges absolutely only when $r>1/2$. We omit the argument as it 
will not affect the subsequent discussion. 

\section{Summary of the  Pseudodifferential Calculus}\label{psido}

\subsection{Basic Definitions and Results}\label{BASIC}
In this Section we introduce the terminology and basic results of the 
Connes-Moscovici pseudodifferential calculus. We will 
describe our version of the asymptotic expansions of this calculus in Section
\ref{higgystricks} and then use them in 
Sections \ref{secresolvent} and \ref{residue}. 
This calculus works in great generality, only needing an unbounded
self-adjoint operator $D$. Just as we did in the remarks following 
Definition \ref{qck}, we set
\ben |D|_1=(1+D^2)^{1/2},\ \ \ \ \delta_1(T)=[|D|_1,T],
\ \ T\in\mbox{dom}\delta.\een
We follow the discussion of the pseudodifferential 
calculus in \cite{C4}, 
using $|D|_1$ and $\delta_1$, instead of $|D|$ and $\delta$. 
In order to ensure that the calculus works in this modified setting
we flesh out explanations in \cite{C4} and record some 
elementary properties which are trivial to prove, 
but are often used without comment. The most important results are
Proposition \ref{L&R} and its Corollary \ref{COR}.
Those familiar with \cite{C4} or \cite{CM} can skip to section 
\ref{higgystricks}.   

So let $D:\mbox{dom}D\subseteq\HH\tto\HH$ be an unbounded self-adjoint 
operator on the Hilbert space $\HH$. For all $k\geq 0$, we set 
\ben \HH_{k}={\rm dom}(1+D^2)^{k/2}={\rm dom}|D|^k\subseteq \HH\een
 and $\HH_\infty=\cap_{k\geq 0}\HH_k$.
Recall that the graph norm topology makes $\HH_{k}$ 
into a Hilbert space with norm $\n\cdot\n_k$ given by
\ben \n\xi\n_k^2=\n\xi\n^2+\n(1+D^2)^{k/2}\xi\n^2\een
where $\n\cdot\n$ is the norm on $\HH$.

We assume that all of our operators $T$, in particular $D$, 
are affiliated to $\cn$ and as in \cite{C4} that 
$T:\HH_\infty\to\HH_\infty$. In this way, all computations 
involving bounded or unbounded operators make sense on the dense 
subspace $\HH_\infty$.

\begin{defn}  For $r\in{\R}$, let 
$op^r$ be the linear space of  operators affiliated to $\cn$ and mapping 
$\HH_\infty\to\HH_\infty$ which are continuous in the norms 
$(\HH_{\infty},\n\cdot\n_k)\to (\HH_{\infty},\n\cdot\n_{k-r})$ for all $k$ 
such that $k-r\geq 0$.
\end{defn}
{\bf Example} The operators $|D|^r$ and $(1+D^2)^{r/2}$ are in $op^r$. 

\begin{lemma}\label{sigmabdd}[compare Lemma 1.1 of \cite{C4}] Let 
$b\in\cap_{n\geq 0}{\rm dom}\delta_1^n$. With $\s_1(b)=|D|_1b|D|_1^{-1}$ 
and $\varepsilon_1(b)=\delta_1(b)|D|_1^{-1}$ we have

1) $\s_1=Id+\varepsilon_1$,

2) $\varepsilon_1^n(b)=\delta_1^n(b)|D|_1^{-n}\in\cn\ \ \ \forall n$,

3) $\s_1^n(b)=(Id+\varepsilon_1)^n(b)=
\sum_{k=0}^n{n\choose k}\delta_1^k(b)
|D|_1^{-k}\in\cn\ \ \ \forall n.$
\end{lemma}
\begin{proof}
The first statement is straightforward. The second follows because $\delta_1$ 
is a derivation with $\delta_1(|D|_1)=0$. The third is just the binomial 
theorem applied to 1).
\end{proof}
Similarly, if $b\in op^0$, $\s_1^{-n}(b):=|D|_1^{-n}b|D|_1^n\in\cn$ for all 
$n$ and
\ben |D|_1^{-n}b|D|_1^n=\sum_{k=0}^n{n\choose k}
|D|_1^{-k}\delta_1^k(b). \een
\begin{cor}\label{smoothiszero} If $b\in\cap_{n\geq 0}{\rm dom}\delta_1^n$ 
then $b\in op^0$.
\end{cor}
\begin{proof}
For $k$ an integer and $\xi\in\HH_\infty$,
\ben \n b\xi\n_k^2=\n b\xi\n^2+\n|D|_1^kb\xi\n^2
=\n b\xi\n^2+\n\s_1^k(b)|D|^k_1\xi\n^2
\leq C(\n\xi\n^2+\n|D|_1^k\xi\n^2).\een
The case of general $k\in{\R}$ now follows by interpolation.
\end{proof}
Observe that by the above Lemma, if $b\in op^0$ then $b-\s_1(b)=
-\varepsilon_1(b)=-\delta_1(b)|D|_1^{-1}\in op^{-1}$. Thus if 
$(\A,\HH,\D)$ is a $QC^\infty$ spectral triple, and $b=a$ or $[\D,a]$ for 
$a\in\A$, then (with $\D$ playing the role of $D$) $b\in op^0$ and 
$b-|\D|_1b|\D|_1^{-1}\in op^{-1}$. 

{\bf Example} In even the most elementary case $\A=C^\infty(S^1)$, 
$\HH=L^2(S^1)$, $a=M_z$, the operator of multiplication by $z$, and $\D=
\frac{1}{i}\frac{d}{d\theta}$ one can easily see that 
$a\in \cap_{n\geq 0}\mbox{dom}\delta_1^n$ but that $[\D^2,a]$ is not bounded. 
In general, $[\D^2,a]$ is about the same size as $|\D|.$ 

\begin{defn} We define the commuting
operators $L_1,R_1$ on the space of operators on $\mathcal H_{\infty}$ by
\ben L_1(T)=(1+D^2)^{-1/2}[D^2,T]=|D|_1^{-1}[|D|_1^2,T],\een
\ben R_1(T)=[D^2,T](1+D^2)^{-1/2}=[|D|_1^2,T]|D|_1^{-1}.\een
\end{defn}

\begin{prop}\label{L&R}[Compare Lemma 2 \cite{C4}] For all 
$b\in op^0$ the following are equivalent:

1) $b\in\bigcap_{n\geq 0}{\rm dom}\delta_1^n$,

2) $b\in\bigcap_{k,l\geq 0}{\rm dom}L_1^k\circ R_1^l$.
\end{prop}

\begin{proof} First observe that $L_1,R_1$ and $\delta_1$ are mutually 
commuting as maps on the space of operators on $\HH_\infty$ which is dense in
$\mathcal H$.
To see that 1) implies 2), let $b\in\cap_{n\geq 0}\mbox{dom}\delta_1^n$ and 
observe $L_1(b)=\delta_1(b)+\s_1^{-1}(\delta_1(b))\in op^0$ by Corollary 
\ref{smoothiszero}. Similarly $R_1(b)=\delta_1(b)+\s_1(\delta_1(b))\in op^0$, 
thus $b\in\mbox{dom}R_1\cap\mbox{dom}L_1$.
Now $\s_1$ and $\delta_1$ commute and leave 
$\cap_{n\geq 0}\mbox{dom}\delta_1^n$ 
invariant, . Since
$L_1^m=\sum_{k=0}^m{m\choose k}\s_1^{-k}\delta_1^m,$
and we easily see that $b$ is in the domain of all powers of $L_1$. A little 
more work shows that $b\in\cap_{k,l\geq 0}\mbox{dom}L^k_1\circ R_1^l$.

As noted in \cite{C4} the implication $2)\Rightarrow 1)$ is more subtle. 
We first want to see that $b\in \cap_{k,l\geq 0}\mbox{dom}L^k_1\circ R_1^l$ 
implies that $b\in\mbox{dom}\delta_1$. Using
\ben |D|_1=\frac{1}{\pi}\int_0^\infty (1+D^2)(1+D^2+u)^{-1}u^{-1/2}du\een

we show by a similar calculation to Lemma 2 of \cite{C4}: 
\bean [|D|_1,b]&=&\frac{1}{\pi}\int_0^\infty 
\left(R_1(b)|D|_1(1+D^2+u)^{-2}u^{1/2} +
 h_u(D)L_1^2(b)(1+D^2+u)^{-1}u^{1/2}\right)du\eean
where $h_u(D)=\left(u(1+D^2+u)^{-1}-1\right)$. We  
show that the integrals of both terms are bounded. For the first, $R_1(b)$ 
is bounded so we compute
$$\frac{1}{\pi}\int_0^\infty|D|_1(1+D^2+u)^{-2}u^{1/2}du=|D|_1
\frac{1}{\pi}\int_0^\infty(1+D^2+u)^{-2}u^{1/2}du$$
$$\leq|D|_1\frac{1}{\pi}\int_0^\infty(1+D^2+u)^{-1}u^{-1/2}du
=|D|_1|D|_1^{-1}=1.$$
Here we used the formula $x^{-1/2}=
\frac{1}{\pi}\int_0^\infty(x+u)^{-1}u^{-1/2}du$. 
To see that the second term is bounded we observe that $\n h_u(D)\n\leq 1$, 
$L^2_1(b)$ is bounded by hypothesis, and the integral $\int_0^\infty
(1+D^2+u)^{-1}u^{1/2}du$ converges absolutely by estimating $\int_0^1$ and 
$\int_1^\infty$ separately.

The conclusion is that $b\in\cap_{l,k}\mbox{dom}L_1^k\circ R_1^l$ 
implies $b\in\mbox{dom}\delta_1$. As 
$L_1(b)\in \cap_{l,k}\mbox{dom}L_1^k\circ R_1^l$, we have 
$L_1(b)\in\mbox{dom}\delta_1$. Since 
$\delta_1(L_1(b))=L_1(\delta_1(b))$, we have $\delta_1(b)\in\mbox{dom}L_1$. 
Similarly $\delta_1(b)\in\cap_{k,l}\mbox{dom}L_1^k\circ R_1^l$ and so by the 
above $\delta_1(b)\in\mbox{dom}\delta_1$; i.e. $b\in\mbox{dom}\delta_1^2$. We 
continue by induction.
\end{proof} 

\begin{defn} For $r\in{\R}$
\ben OP^r=|D|_1^r\left(\bigcap_{n\geq 0}{\rm dom}\delta_1^n\right)\subseteq 
op^r\cdot op^0\subseteq op^r.\een
If $T\in OP^r$ we say that the order of $T$ is (at most) $r$.
The definition is actually symmetric, since for $r$ an integer 
(at least) we have by Lemma \ref{sigmabdd}
\bean OP^r=|D|_1^r\left(\bigcap\mbox{dom}\delta_1^n\right)
=|D|_1^r\left(\bigcap\mbox{dom}\delta_1^n\right)
|D|_1^{-r}|D|_1^r
\subseteq\left(\bigcap\mbox{dom}\delta_1^n\right)|D|_1^r.\eean
From this we easily see that $OP^r\cdot OP^s\subseteq OP^{r+s}.$
Finally, we note that if $b\in OP^r$ for $r\geq 0$, then since $b=|D|_1^ra$ 
for some $a\in OP^0$, we get
$[|D|_1,b]=|D|_1^r[|D|_1,a]=|D|_1^r\delta_1(a),$
so $[|D|_1,b]\in OP^r$.
\end{defn}
{\bf Remarks} An operator $T\in OP^r$ if and only if 
$|D|_1^{-r}T\in\cap_{n\geq 0}\mbox{dom}\delta_1^n$. 
 Observe that operators of order at most zero are 
bounded. If $|D|_1^{-1}$ is $p$-summable and $T$ has order $-n$ then, $T$
is $p/n$-summable.

{\bf Important Observations}\\ 
1) If $f$ is a bounded Borel function then $f(D)\in\cn$ and 
$\delta_1(f(D))=0$, implies $f(D)\in OP^0.$\\
2) If $g$ is an unbounded Borel function such that $1/g$ is 
bounded on $spec(D)$ and both $g(D)|D|_1^{-1}$ and $g(D)^{-1}|D|_1$ 
are bounded, then for each $r$, $OP^r=|g(D)|^rOP^0$. This follows since 
$OP^0$ is an algebra and both $|g(D)|^r|D|_1^{-r}$ and 
$g(D)^{-r}|D|_1^r$ are in $OP^0$. 
We note that if $|D|$ is {\bf not} invertible then we get strict containment
$|D|^r OP^0\subset|D|_1^r OP^0$.
These observations prove the next Lemma.

\begin{lemma}\label{OPs} If $\mu\in{\C}$ is in the 
resolvent set of $D^2$ then 
\ben OP^r=|(\mu-D^2)^{1/2}|^r\left(\bigcap_{n\geq 0}{\rm dom}\delta_1^n\right)
.\een
\end{lemma}

\begin{cor}\label{COR}\label{expect} Let $(\A,\HH,\D)$ be 
a $QC^\infty$ spectral triple, and 
suppose $a\in\A$. Then for $n\geq 0$, $a^{(n)}$ and $[\D, a]^{(n)}$ are 
in $OP^{n}$.
\end{cor}
\begin{proof}
Writing $T^{(n)}=[D^2,[D^2,[\cdots[D^2,T]\cdots]]]$ 
for the $n$-th iterated commutator with $D^2$,
we have already observed that 
$L_1^2(b)=|D|_1^{-2}[D^2,[D^2,b]]=|D|_1^{-2}b^{(2)}.$
Similarly, $L_1^n(b)=|D|_1^{-n}b^{(n)}$. So if 
$b\in\cap_{k\geq 0}\mbox{dom}\delta_1^k$, then $L_1^n(b)$ is also in 
$\cap_{k\geq 0}\mbox{dom}\delta_1^k$ and hence $b^{(n)}=|D|_1^nL_1^n(b)\in 
OP^n$. 
\end{proof}

\subsection{The Pseudodifferential Expansion}\label{higgystricks}

Next we describe the asymptotic expansions introduced by 
Connes and Moscovici in 
\cite{C4,CM}. Their principal result is that if $T\in OP^k$ for $k$ integral, 
then 
for any $z\in{\C}$ 
\bean 
(1+D^2)^{z}T&=&T(1+D^2)^{z}+zT^{(1)}(1+D^2)^{z-1}+
\frac{z(z-1)}{2}T^{(2)}(1+D^2)^{z-2}+\cdots\nno
&&\cdots + \frac{z(z-1)\cdots(z-n+1)}{n!}T^{(n)}(1+D^2)^{z-n}+P,\eean
where $P\in OP^{k-(n+1)+2Re(z)}$. This result is proved in both of the papers 
\cite{C4,CM}, but subsequently a simpler proof has been given by Higson
\cite{hig}.

Because of Higson's idea, we will not need the full force of this expansion, 
but only a simple algebraic version of it (Lemmas \ref{higsexpan} and
\ref{firstexpan}). We briefly 
sketch the idea behind Higson's proof in the way in which we will 
use it. So we suppose that we have a $QC^\infty$ (odd) spectral triple 
$(\A,\HH,\D)$ with dimension $p\geq 1$. We use $\D$ to define the 
pseudodifferential calculus on $\HH$ as in the previous section. 
Let $Q=(1+s^2+D^2)$ where $D$ is $\tD$ as defined in Section 5.3
and where $s\in [0,\infty)$.
For $Re(z)>p/2$ we write $Q^{-z}$ using Cauchy's formula
\ben Q^{-z}=\frac{1}{2\pi i}\int_l\lambda^{-z}(\lambda-Q)^{-1}d\lambda,\een
where $l$ is a vertical line $\lambda=a+iv$ parametrized by $v\in\R$ 
with $0<a<1/2$ fixed. One checks that the integral indeed converges in operator 
norm (to an element of $\cn$) and by using the spectral theorem for $Q$
(in terms of its spectral resolution) it converges to $Q^{-z}$ (principal 
branch). Computing commutators of $Q^{-z}$ with an operator $T\in OP^m$ then 
reduces to an iterative calculation of commutators with $(\lambda-Q)^{-1}$. 
The exact result we need is the following.

\begin{lemma} \label{higsexpan} Let $m,n,k$ be non-negative integers and 
$T\in OP^m$. Then 
\bean (\lambda-Q)^{-n}T&=&T(\lambda-Q)^{-n}+nT^{(1)}(\lambda-Q)^{-(n+1)}+
\frac{n(n+1)}{2}T^{(2)}(\lambda-Q)^{-(n+2)}+\cdots\nno
&\cdots&+{n+k-1\choose k}T^{(k)}
(\lambda-Q)^{-(n+k)}+P(\lambda)\nno
&=& \sum_{j=0}^{k} {n+j-1\choose j}T^{(j)}
(\lambda-Q)^{-(n+j)} + P(\lambda)\eean
where the remainder $P(\lambda)$ has order $-(2n+k-m+1)$ and is given by
\ben P(\lambda)=\sum_{j=1}^n {j+k-1\choose k}
(\lambda-Q)^{j-n-1}T^{(k+1)}(\lambda-Q)^{-j-k}.\een
\end{lemma}
\begin{proof} The proof is inductive using the resolvent formula
\ben (\lambda-Q)^{-1}T=T(\lambda-Q)^{-1}+(\lambda-Q)^{-1}T^{(1)}
(\lambda-Q)^{-1}.\een
For example with $k=0$
\bean(\lambda-Q)^{-n}T&=&(\lambda-Q)^{-(n-1)}\left(T(\lambda-Q)^{-1}+
(\lambda-Q)^{-1}T^{(1)}(\lambda-Q)^{-1}\right)\nno
&=&\vdots\nno
&=&T(\lambda-Q)^{-n}+\left[(\lambda-Q)^{-1}T^{(1)}(\lambda-Q)^{-n}+\cdots\right.
\nno
&\cdots&\left.+(\lambda-Q)^{-(n-1)}T^{(1)}(\lambda-Q)^{-2}+
(\lambda-Q)^{-n}T^{(1)}(\lambda-Q)^{-1}\right],\eean
the term in square brackets being the remainder, which clearly has order 
$-2(n+1)+m+1=-(2n+0-m+1)$.
For $k=1$, we move $T^{(1)}$ to the left in each term of the remainder using 
the resolvent formula again. We arrive at 
\ben (\lambda-Q)^{-n}T=T(\lambda-Q)^{-n}+nT^{(1)}(\lambda-Q)^{-(n+1)}+
\sum_{j=1}^nj(\lambda-Q)^{-(n+1-j)}T^{(2)}(\lambda-Q)^{-(j+1)},\een
and the new remainder has order $(m+2)-2(n+2)=-(2n+1-m+1)$. While it is clear 
how to proceed, in order to keep track of the constants in the remainder terms
one must use the formula:
$$\sum_{j=1}^n {j+k-1\choose k}=
{n+k\choose k+1}$$
which is easily proved by induction on $n$ with $k$ held fixed. 
\end{proof}

\begin{cor*}
Let $n,M$ be positive integers and $A\in OP^k$. Let $R=(\lambda-Q)^{-1}$.
Then,
$$R^n AR^{-n} = \sum_{j=0}^{M} 
{n+j-1\choose j}A^{(j)}R^j  + P$$
where
$$P=\sum_{j=1}^n {j+M-1\choose M}R^{n+1-j}
A^{(M+1)}R^{M+j-n}$$
and $P$ has order $k-M-1$.
\end{cor*}

\begin{lemma}\label{normtrick} Let $k,n$ be non-negative integers, $s\geq 0$, 
and suppose
$\lambda\in{\C}$, $a=Re(\lambda)$ with $0<a<1/2$. Then for $A \in OP^k$ and with
$R_s(\lambda)=(\lambda-(1+D^2+s^2))^{-1}$ we have
$$ \n R_s(\lambda)^{n/2+k/2}AR_s(\lambda)^{-n/2}\n\leq C_{n,k}
\;\;and\;\;\n R_s(\lambda)^{-n/2}AR_s(\lambda)^{n/2+k/2}\n\leq C_{n,k}$$
where $C_{n,k}$ is constant independent of $s$ and $\lambda$ (square roots use 
the principal branch of $\log$.)
\end{lemma}
\begin{proof} Since $(OP^k)^* = OP^k$ and $(R_s(\lambda))^* = R_s(\bar\lambda)$
we only need to prove the first inequality. We begin by using the numerical 
inequality $(a^2+b^2)^{1/2}\leq |a|+|b|$ repeatedly to see that
\ben |R_s(\lambda)|^{-1/2}=(v^2+((1-a)+D^2+s^2)^2)^{1/4}\leq\cdots
\leq(|v|^{1/2}+(1-a)^{1/2}+|D|+s).\een
We then use the companion inequality $|a|+|b|\leq\sqrt{2}(a^2+b^2)^{1/2}$ and 
the fact that $1<\sqrt{2}$ repeatedly to see that
\ben (|v|^{1/2}+(1-a)^{1/2}+|D|+s)
\leq\cdots\leq 4(v^2+((1-a)+D^2+s^2)^2)^{1/4}=4|R_s(\lambda)|^{-1/2}.\een
Hence, as non-negative bounded operators in $OP^0$ we have for $j\geq 0,$
\ben (|v|^{1/2}+(1-a)^{1/2}+|D|+s)^{-j}|R_s(\lambda)|^{-j/2}\leq 1\een
\ben |R_s(\lambda)|^{j/2}((1-a)^{1/2}+|D|)^j\leq|R_s(\lambda)|^{j/2}
(|v|^{1/2}+(1-a)^{1/2}+|D|+s)^j\leq 4^j.\een
Now, $A=((1-a)^{1/2}+|D|)^k B$ for $B\in OP^0$ by Lemma \ref{OPs}. 
Letting $|D|_a=(1-a)^{1/2}+|D|$, making the substitution $A=|D|_a^k B$,
and using the identity (valid for any constant, $c$):
$$B(c+|D|)^n=\sum_{j=0}^{n}{n \choose j}
(c+|D|)^j(-\delta)^{n-j}(B)$$
we have,
\bean &&\n R_s(\lambda)^{k/2+n/2}AR_s(\lambda)^{-n/2}\n =
\n |D|_a^kR_s(\lambda)^{k/2+n/2}BR_s(\lambda)^{-n/2}\n\nno
&\leq& \n |D|_a^kR_s(\lambda)^{k/2+n/2}
B(|v|^{1/2}+(1-a)^{1/2}+|D|+s)^n\n\nno
&=&\n |D|_a^kR_s(\lambda)^{k/2+n/2}\sum_{j=0}^{n}{n \choose j}
(|v|^{1/2}+(1-a)^{1/2}+|D|+s)^{j}(-\delta)^{n-j}(B)\n\nno
&\leq& \sum_{j=0}^{n}{n\choose j}\n |D|_a^k R_s(\lambda)^{k/2}\n\cdot 
\n R_s(\lambda)^{n/2}
(|v|^{1/2}+(1-a)^{1/2}+|D|+s)^j\n\cdot \n\delta^{n-j}(B)\n\nno
&\leq& \sum_{j=0}^{n}{n\choose j} 4^k \n  
R_s(\lambda)^{n/2-j/2}
\n\cdot \n R_s(\lambda)^{j/2}(|v|^{1/2}+(1-a)^{1/2}+|D|+s)^j\n\cdot
\n\delta^{n-j}(B)\n\nno
&\leq& \sum_{j=0}^{n}{n\choose j}
\left(\frac{1}{1-a}\right)^
{(n-j)/2} 4^{k+j} \n\delta^{n-j}(B)\n.\eean
\end{proof}

\begin{rems*}
We are thinking of the operator $X$ in the following lemmas as $\{\tD,q\}$
in Section \ref{G}
so that indeed $s^2+s\{\tD,q\}+\tD^2 = (\tD + sq)^2\geq 0$. With this 
hypothesis on $X$ and with $0<a=Re(\lambda)<1/2$ we see that the spectrum
of $(1+s^2+sX+\tD^2)$ is bounded away
from the line $Re(\lambda)=a$ by at least $(1-a)$ independent of $s$ and 
$\lambda$. This hypothesis is crucial in Lemma \ref{uni}. Lemma
\ref{higsexpan} and the next lemma form  
the algebraic heart of our pseudodifferential expansion.
\end{rems*}

\begin{lemma}\label{firstexpan} Let $A_i\in OP^{n_i}$ for $i=1,...,m$ and let
$0<a=Re(\lambda)<1/2$ as above. We consider the operator 
$$R_s(\lambda)A_1R_s(\lambda)A_2R_s(\lambda)\cdots R_s(\lambda)A_m
\tilde R_s(\lambda)$$
where $R_s(\lambda)=(\lambda-(1+s^2+D^2))^{-1}$ and
$\tilde R_s(\lambda)=(\lambda-(1+s^2+sX+D^2))^{-1}$
where $X$ is self-adjoint, bounded and $s^2+sX+D^2\geq 0$. 
Then for all $M\geq 0$
\ben R_s(\lambda)A_1R_s(\lambda)A_2\cdots 
A_m\tilde R_s(\lambda)=\sum_{|k|=0}^MC(k)A_1^{(k_1)}\cdots 
A_m^{(k_m)}R_s(\lambda)^{m+|k|}\tilde R_s(\lambda)+P_{M,m}\een
where $P_{M,m}$ is of order (at most) $-2m-M-3+|n|$, and  $k$ and $n$
are multi-indices with $|k|=k_1+\cdots+k_m$ and $|n|= n_1+\cdots+n_m$.
The constant $C(k)$ is given by
\ben C(k)=\frac{(|k|+m)!}{k_1!k_2!\cdots k_m!(k_1+1)(k_1+k_2+2)\cdots(|k|+m)}=
(|k|+m)!\alpha(k).
\een
\end{lemma}
\begin{proof} In order to lighten the notation, we abbreviate $R:=R_s(\lambda)$
and $\tilde R:=\tilde R_s(\lambda).$
Then we write:
\be 
\label{res}
RA_1RA_2\cdots A_m\tilde R=RA_1R^{-1}R^2 A_2R^{-2}\cdots R^mA_mR^{-m}R^m
\tilde R\ee
so that we have a product of $m$ operators $R^jA_jR^{-j}$ and an additional 
factor of $R^m\tilde R$. 

Now by the corollary to Lemma \ref{higsexpan}, if $i$ is a positive integer
$$R^iA_iR^{-i}=\sum_{j=0}^M 
{i+j-1\choose j}A_i^{(j)}R^j +P$$
where $P$ is order $n_i-M-1$.
So if we expand the term $RA_1R^{-1}$ on the right hand side  of 
Equation (\ref{res})
we obtain
\bean &&RA_1R^{-1}R^2A_2R^{-2}\cdots
R^mA_mR^{-m}R^{m}\tilde R\nno
&=&\sum_{k_1=0}^{M}A_1^{(k_1)}R^{2+k_1}A_2R^{-2}\cdots 
A_mR^{-m}R^{m}\tilde R+P_1\nno
&=&\sum_{k_1=0}^{M}A_1^{(k_1)}R^{2+k_1}A_2R^{-2-k_1}\cdots 
R^{m+k_1}A_mR^{-m-k_1}R^{m+k_1}\tilde R+P_1.\eean
Here $P_{1}$ is of order $|n|-2m-M-3$. In fact, 
$$P_1 = RA_1^{(M+1)}R^{M+1}RA_2R\cdots RA_m\tilde R.$$
 Now, expanding this new second term, $R^{2+k_{1}}A_2R^{-(2+k_1)}$, gives
$$\sum_{k_1=0}^{M}\sum_{k_2=0}^{M}C_1(k)
A_1^{(k_1)}A_2^{(k_2)}R^{3+k_1+k_2}A_3R^{-3-k_1-k_2}
\ldots R^{-m-k_1-k_2}R^{m+k_1+k_2}\tilde R +P_1+P_2$$
where $C_1(k)={1+k_1+k_2\choose k_2}$ and $P_2$ is of order at most 
$$\mbox{max}(|n|-M-2m-3-k_1)=|n|-M-2m-3.$$
In fact,
$$P_2=\sum_{k_1=0}^M \sum_{j=1}^{2+k_1}{j+M-1\choose M} A_1^{(k_1)}
R^{3+k_1 -j}A_2^{(M+1)}R^{M+j}RA_3R\cdots RA_m\tilde R.$$
Repeating this argument shows that
\ben
RA_1RA_2\cdots A_m\tilde R=\sum_{k_1,k_2,\ldots,k_m=0}^MC(k)A_1^{(k_1)}\cdots 
A_m^{(k_m)}R^{m+|k|}\tilde R+ \mbox{ order }(|n|-M-2m-3)\een
where 
$$C(k)={k_1\choose k_1}{1+k_1+k_2\choose k_2}\ldots
{m-1+k_1+\ldots+k_m\choose k_m}.$$
What remains is to eliminate excess terms in the previous sum
and determine the coefficient $C(k)$.
Observe first that if $k_1+\ldots+k_m>M$ then the order of
$A_1^{(k_1)}\cdots A_m^{(k_m)}R^{m+|k|}\tilde R$ is $|n|-2-2m-|k|<|n|-2-2m-M$
which in turn is $\leq |n|-3-2m-M$ (since only integers are considered here.) 
Thus
\ben
RA_1RA_2\cdots A_m\tilde R=\sum_{|k|=0}^M C(k)A_1^{(k_1)}\cdots 
A_m^{(k_m)}R^{m+|k|}\tilde R +P_{M,m} \een
where $P_{M,m}$ has order $(|n|-M-2m-3).$ 
Finally $C(k)$ is 
\bean 
&&\frac{(1+k_1+k_2)!}{k_2!(1+k_1)!}\frac{(2+k_1+k_2+k_3)!}{k_3!(2+k_1+k_2)!}
\cdots 
\frac{(m-1+k_1+\cdots+k_m)!}{k_m!(m-1+k_1+\cdots+k_{m-1})!}\nno
&&\qquad\quad=\frac{k_1!}{k_1!}\frac{(1+k_1+k_2)!}{k_2!(1+k_1)!}
\frac{(2+k_1+k_2+k
_3)!}{k_3!(2+k_1+k_2)!}\cdots 
\frac{(m-1+k_1+\cdots+k_m)!}{k_m!(m-1+k_1+\cdots+k_{m-1})!}\nno
&&\qquad\quad= \frac{(m-1+k_1+\cdots+k_{m})!}{k_1!\cdots 
k_m!(1+k_1)(2+k_1+k_2)\cdots(m-1+k_1+\cdots+k_{m-1})}\nno
&&\qquad\quad=\frac{(m+|k|)!}{k_1!\cdots 
k_m!(1+k_1)(2+k_1+k_2)\cdots(m+k_1+\cdots+k_{m})}.\eean
\end{proof}

\begin{lemma}\label{uni} With the assumptions and notation of the last Lemma
including the assumption that $A_i\in OP^{n_i}$ for each $i$, 
there is a positive constant $C$ such that
\ben \n (\lambda-(1+D^2+s^2))^{m+M/2+1/2-|n|/2}P_{M,m}\n\leq C\een
independent of $s$ and $\lambda$ (though it depends on $M$ and $m$ and 
the $A_i$). If the final $\tilde{R}_s(\lambda)$ is actually an $R_s(\lambda)$
then we can replace the $1/2$ in the exponent with $3/2$: this is important
in the proof of Proposition 8.1.
\end{lemma}
\begin{proof} 
The remainder $P_{M,m}$ in the previous lemma
obtained after applying the pseudodifferential expansion
has terms of two kinds. The first kind we consider are the bookkeeping terms
at the end of the proof of the last lemma. They are of the form
$$P=A_1^{(k_1)}\cdots A_m^{(k_m)}R^{m+|k|}\tilde R$$
with $|k|>M$. Then $R^{-(m+M/2+1/2)+|n|/2}P$ is uniformly bounded
by an application of Lemma \ref{normtrick}.

The other terms are the ones $P_1, P_2,\cdots,P_m$ obtained in the proof of 
the last lemma. Recall:
$$P_1 = RA_1^{(M+1)}R^{M+1}RA_2R\cdots RA_m\tilde R$$
while a typical summand of $P_2$ is an integer multiple of:
$$A_1^{(k_1)}
R^{3+k_1 -j}A_2^{(M+1)}R^{M+j}RA_3R\cdots RA_m\tilde R$$
where $1\leq j\leq 2+k_1$ and $0\leq k_1\leq M$.

Similarly, for $1\leq i\leq m$, a typical summand of $P_i$ is an integer 
multiple of:
$$A_1^{(k_1)}A_2^{(k_2)}\cdots A_{i-1}^{(k_{i-1})}
R^{i+1+k_1+k_2+\cdots k_{i-1}-j}A_i^{(M+1)}
R^{M+j}RA_{i+1}R\cdots RA_m\tilde R$$
where $1\leq j\leq i+k_1+k_2+\cdots+k_{i-1}$ and 
$0\leq k_1,k_2,\cdots, k_{i-1} \leq M.$

We work with the typical summand of $P_i$ above, and
let $B=A_1^{(k_1)}A_2^{(k_2)}\cdots A_{i-1}^{(k_{i-1})}$ which has order 
$(k_1+k_2+\cdots k_{i-1})+(n_1+n_2+\cdots n_{i-1})=|k|+|n|_{i-1},$
where we have used the notation $|n|_j=n_1+n_2+\cdots n_j$. We will also
use the notation $|n|^{j+1}=n_{j+1}+\cdots n_m$ so that $|n|=|n|_j+|n|^{j+1}$.
We need to show that:
$$R^{-(m+M/2+1/2)+|n|/2}BR^{i+1+|k|-j}A_i^{(M+1)}R^{M+j}RA_{i+1}R
\cdots RA_m\tilde R$$
is bounded independent of $s$ and $\lambda$. 
So, we calculate:
\bean &&R^{-(m+M/2+1/2)+|n|/2}BR^{i+1+|k|-j}
A_i^{(M+1)}R^{M+j}RA_{i+1}R\cdots RA_m\tilde R\\
&=&R^{-(m+M/2+1/2-|n|/2)}BR^{(|k|+|n|_{i-1})/2}R^{(m+M/2+1/2-|n|/2)}
R^{|k|/2}R^{-(m+M/2+j-i-1/2)+|n|^{i}/2}\\
&&\times A_i^{(M+1)}R^{M+j}RA_{i+1}R\cdots RA_m\tilde R\\
&=&\left(R^{-(m+M/2+1/2-|n|/2)}BR^{(|k|+|n|_{i-1})/2}R^{(m+M/2+1/2-|n|/2)}
\right)R^{|k|/2}\\
&&\times\left(R^{-(M/2+m+j-i-1/2)+|n|^{i}/2}A_i^{(M+1)}
R^{[(M/2+m+j-i-1/2)-|n|^{i}/2]+(n_i+M+1)/2}\right)\\
&&\times (R^{-(m-i-1)+|n|^{i+1}/2}A_{i+1}R^{(m-i-1-|n|^{i+1}/2)+n_{i+1}/2})
\cdots\times\\
&&\times (R^{-1+(n_{m-1}+n_m)/2}A_{m-1}R^{1-(n_{m-1}+n_m)/2+n_{m-1}/2})
(R^{n_m/2}A_m)\tilde R.\eean
Then each bracketed term in the last expression is bounded independent of
$s$ and $\lambda$ by an application of Lemma \ref{normtrick}, while
$$\n R^{|k|/2}\n \leq \left(\frac{1}{1-a}\right)^{|k|/2}\;\;and\;\;
\n \tilde R\n \leq \frac{1}{1-a}$$
by Lemma \ref{normlambda} parts (a) and (b) and the condition on the 
operator $B$.
We remark that the cases when $|k|=0$ are included in this calculation.
\end{proof}
\section{The Resolvent Cocycle}\label{secresolvent}

The standing assumptions for the computations in this Section and Section 
\ref{residue} are that 
$(\A,\HH,\D)$ is a $QC^\infty$  
semifinite spectral triple.  
We denote  the dimension of $(\A,\HH,\D)$ by $p$, 
and we suppose that $p\geq 1$.
For $r>0$, $(1+\D^2)^{-p/2-r}$ is trace class, and if $u\in\A$ is unitary, 
the 
spectral flow from $\D$ to $u\D u^*$ is given by
\ben 
\frac{1}{C_{p/2+r}}\int_0^1\tau\left(u[\D,u^*](1+(\D+tu[\D,u^*])^2)^{-p/2-r})
\right)dt,\een
where the constant is 
\ben 
C_{p/2+r}=\int_{-\infty}^\infty(1+x^2)^{-p/2-r}dx=\frac{\Gamma((p-1)/2+r)
\Gamma(1/
2)}{\Gamma(p/2+r)}.\een
We saw at the end of Section 5 that 
 we may replace this formula with
\be C_{p/2+r}sf(\D,u^*\D u)=\int_0^\infty 
S\tau(q(1+\tD^2+s^2+s\{\tD,q\})^{-p/2-r})ds,\label{ancont}\ee
for any $r>0$
where  we remind 
the reader that we are now computing $sf(\D,u^*\D u)=-sf(\D,uD 
u^*)$. {\em It is important to observe that the left hand side of equation 
(\ref{ancont}) provides a meromorphic continuation of the function of $r$ 
defined by the integral on the right hand side.}

This meromorphic continuation  has simple poles at $r=(1-p)/2-k$, $k=0,1,2,...$.
 The 
residue of the left hand side at $r=(1-p)/2$ is precisely $sf(\D,u^*\D u)$ so
 we may write
\bean sf(\D,u^*\D u)&=& \frac{1}{2\pi i}\int_\gamma C_{p/2+z}sf(\D,u^*\D u)dz
\eean
where $\gamma=\{z=(1-p)/2+\epsilon e^{i\theta},\ 0\leq\theta\leq 2\pi\}$ for 
a 
suitably small $\epsilon$.  
Our aim now is  to compute these residues from the integral formula.
We note however that our previous estimates (Sections \ref{nandt}, \ref{tande}) 
do not 
allow us to use the formula in Equation (\ref{ancont}) for the analytic 
continuation of the right hand side to a deleted neighbourhood of the 
critical point  $r=(1-p)/2$.

In order to obtain such a formula we need to perform further manipulations
on this integral.
These involve expanding the integrand 
using the Cauchy formula, the resolvent expansion, 
and the pseudodifferential calculus. 
Then we make some estimates to show that we can discard remainders in these 
expansions and that the terms left over do indeed analytically 
continue to a deleted neighbourhood of the critical point $r=(1-p)/2$. These 
continuations exist by virtue of the isolated spectral dimension hypothesis.
We note some facts here.

{\bf Observation 1}.
Let $N=[p/2]+1$ be the least integer strictly greater than  $p/2\geq1/2$. 
Observe that if $p$ is 
an odd integer, $N=(p+1)/2$, while if $p$ is an even integer, $N=p/2+1$.
More generally there is always a $\delta>0$ so that $N\geq (p+\delta)/2$.
 
{\bf Observation 2}.
Let $l$ be the vertical line $\{\lambda:Re(\lambda)=a\}$ where $0<a<1/2$ is 
fixed. We have already observed in a previous remark that the line $l$ is
uniformly bounded away from the spectrum of $(1+\tD^2+s\{\tD,q\}+s^2)$ by 
$(1-a)$ as $s$ varies since $(1+\tD^2+s\{\tD,q\}+s^2)= 1+(\tD+sq)^2$. This 
implies, in particular, that $\n\tilde R_s(\lambda)\n\leq 
\sqrt{v^2+(1-a)^2}^{-1}\leq 1/(1-a)$ for all $s$ and $\lambda$.
In some of our estimates it also helps to have $\n \{\tD,q\}\n<\sqrt{2}$.
This can be achieved by scaling $\tD$; since both our spectral flow
formula and our final cocycle formulas are invariant under this multiplicative
change of scale we are free to rescale for each individual $q$.

{\bf Observation 3} In Section \ref{G} we began with a $QC^\infty$ 
semifinite spectral triple $(\A,\HH,\D)$ and looked at various operators
affiliated with $\tilde{\mathcal N}=M_2\otimes M_2\otimes \mathcal N$. 
In particular, we considered
$$\tD=\sigma_2\otimes 1_2\otimes \D\;\;,\;\; 
q=\sigma_3\otimes\left(\begin{array}{cc}
                   0 & -iu^{-1} \\
                   iu & 0
                   \end{array} \right)\;\;and\;\;
\{\tD,q\}=\sigma_1\otimes\left(\begin{array}{cc}
                   0 & [\D,u^{-1}] \\
                   -[\D,u] & 0
                   \end{array} \right).$$
Now clearly, $|\tD|_1=1_2\otimes 1_2\otimes |\D|_1$ while both $q$ 
and $\{\tD,q\}$
are in the algebra $M_2\otimes M_2\otimes (OP_{\D}^0)$. Letting 
$\tilde\delta_1$ and $\delta_1$ be the derivations induced by $|\tD|_1$ and
$|\D|_1$, we see that 
$\tilde\delta_1=Id\otimes Id\otimes \delta_1$, so
$$OP_{\tD}^0=M_2\otimes M_2\otimes (OP_{\D}^0).$$
That is, both $q$ and $\{\tD,q\}$ are in $OP^0$ relative to the operator
$\tD$ which satisfies the same summability condition as $\D$.

\subsection{Resolvent expansion of the spectral flow}\label{resexpand} 

We require two estimates to guarantee that various operators which arise from 
the Cauchy formula and the resolvent expansion are trace class. We present 
these as separate lemmas as we will use them repeatedly. The techniques we use 
in these Lemmas are indicative of the methods employed in the remainder of the 
proof.

\begin{lemma}\label{yucko} Let $(\A,\HH,\D)$ be a $QC^\infty$ semifinite 
spectral triple of dimension $p\geq 1$. Let $m$ be a non-negative integer, and 
for $j=0,...,m$ let $A_j\in OP^{0}$. Let $l$ be the vertical line 
$v\mapsto \lambda=a+iv$ for $v\in\R$ and $0<a<1/2$, 
$R_s(\lambda)=(\lambda-(1+s^2+\tD^2))^{-1}$ and 
$\tilde{R}_s(\lambda)=(\lambda-(1+s^2+\tD^2+s\{\tD,q\}))^{-1}$. Then for 
$r\in\bf{C}$ and $Re(r)>0$  
the operator
\ben B(s)=\frac{1}{2\pi i}\int_l\lambda^{-p/2-r}A_0R_s(\lambda)A_1
R_s(\lambda)A_2\cdots R_s(\lambda) A_m\tilde{R}_s(\lambda)d\lambda,\een
is trace class for $m>p/2$ and the function $s^m\n B(s)\n_1$ is integrable on 
$[0,\infty)$ when 
\ben p+\epsilon<1+m\ \ \ \mbox{and}\ \ \ 1+\epsilon<m+2Re(r).\een
\end{lemma}

\begin{proof} By Lemma \ref{firstexpan} we have
\bean B(s)&=&\frac{1}{2\pi i}\int_l\lambda^{-p/2-r}A_0R_s(\lambda)A_1\cdots 
R_s(\lambda)A_m\tilde{R}_s(\lambda)d\lambda\nno
&=&\frac{1}{2\pi i}\int_l\lambda^{-p/2-r}\sum_{|n|=0}^MC(n)A_0A_1^{(n_1)}\cdots 
A_m^{(n_m)}
R_s(\lambda)^{m+|n|}\tilde{R}_s(\lambda)d\lambda\nno
&+& \frac{1}{2\pi i}\int_l\lambda^{-p/2-r}P\tilde{R}_s(\lambda)d\lambda,\eean
where $P$ is of order $-M-1-2m$ and $|n|=n_1+\cdots+n_m$. Now by Lemma
\ref{normtrick},
\ben A_0A_1^{(n_1)}\cdots A_m^{(n_m)}R_s(\lambda)^{|n|/2}\een
is  bounded independently of $s,\lambda$. Thus, provided that $m>p/2$ 
(so that $R_s(\lambda)^{m+|n|/2}$ is trace class); 
using Lemma \ref{uni} to estimate the remainder term; and noting that
since $r$ is complex and $\lambda\in l$,
$$|\lambda^{-p/2-r}|=e^{Re[(-p/2-r)Log(\lambda)]}\leq |\lambda|^{-p/2-Re(r)}
e^{|Im(r)|\pi/2}=C_r|\lambda|^{-p/2-Re(r)},$$
we have
\bea \n B(s)\n_1&\leq& \sum_{|n|=0}^M\frac{C'(n)}{2\pi}\int_{-\infty}^\infty
\sqrt{a^2+v^2}^{-p/2-Re(r)}\n R_s(\lambda)^{m+|n|/2}\n_1
\n\tilde{R}_s(\lambda)\n dv
\nno
&+& \frac{C}{2\pi}\int_{-\infty}^\infty\sqrt{a^2+v^2}^{-p/2-Re(r)}
\n R_s(\lambda)^{m+M/2+1/2}\n_1\n\tilde{R}_s(\lambda)\n dv\nno
&\leq& \sum_{|n|=0}^MC''(n)\left(\int_{-\infty}^\infty
\sqrt{a^2+v^2}^{-p/2-Re(r)}\right.\times\nno
&&\quad\times\left.\sqrt{(1/2+s^2-a)^2+v^2}^{-m-|n|/2+(p+\epsilon)/2}
\sqrt{(1+s^2-a-sc)^2+v^2}^{-1}dv\right)\nno
&+& C'\int_{-\infty}^\infty\left(\sqrt{a^2+v^2}^{-p/2-Re(r)}\right.\times\nno
&&\times\quad\left.\sqrt{(1/2+s^2-a)^2+v^2}^{-M/2-1/2-m+(p+\epsilon)/2}
\sqrt{(1+s^2-a-sc)^2+v^2}^{-1}dv\right),\label{normbee}\eea
The integrals over $v$ are finite for $p/2+Re(r)+m+|n|/2>(p+\epsilon)/2$, which 
(with $r,p,m$ as above) is always true, and $p/2+Re(r)+m+M/2+1/2>(p+\epsilon)/2$, 
which is also always true. 
Now multiply Equation (\ref{normbee}) by $s^m$ and integrate over $[0,\infty)$. 
An application of Lemma \ref{intest} shows that the integral arising from the 
remainder converges if 
\ben M>-m-2+p+\epsilon\ \ \ \mbox{and}\ \ \ M>-m-2Re(r)+\epsilon.\een
Since $M$ is positive, the second holds trivially while the first will hold if
$M\geq p/2-2+\epsilon.$
The other terms are all finite provided
\ben m+1>p+\epsilon\ \ \ \mbox{and}\ \ \ m+2Re(r)>1+\epsilon.\een
This completes the proof.
\end{proof}

\begin{lemma}\label{spaz} Let $(\A,\HH,\D)$ be a $QC^\infty$ semifinite 
spectral triple of dimension $p\geq 1$. Let $m$ be a non-negative integer, and 
for $j=0,...,m$ let $A_j\in OP^{k_j}$, $k_j\geq 0$. Let $l$ be the vertical 
line described above, $R_s(\lambda)=(\lambda-(1+s^2+\tD^2))^{-1}$. 
Then  for $Re(r)>0$ the operator
\ben B(s)=\frac{1}{2\pi i}\int_l\lambda^{-p/2-r}A_0R_s(\lambda)A_1R_s(\lambda)
A_2\cdots R_s(\lambda) A_mR_s(\lambda)d\lambda,\een
is trace class for $Re(r)+m-|k|/2>0$ and the function  
$s^\alpha\times\n B(s)\n_1$ is integrable on $[0,\infty)$ when 
\ben 1+\alpha+|k|-2m<2(Re(r)-\epsilon).\een
\end{lemma}

{\bf Remark} Observe that the estimate in Lemma \ref{spaz} is vastly superior 
to that in Lemma \ref{yucko}. This is because the resolvents are all of the 
same kind, and we may employ Cauchy's formula to do integrals, whereas in Lemma 
\ref{yucko} we are forced to employ trace estimates under the integral. 

\begin{proof} We employ  Lemma \ref{firstexpan} again to find
\bean B(s)&=&\frac{1}{2\pi i}\int_l\lambda^{-p/2-r}A_0R_s(\lambda)A_1
R_s(\lambda)A_2\cdots R_s(\lambda) A_mR_s(\lambda)d\lambda \nno
&=&\frac{1}{2\pi i}\int_l\lambda^{-p/2-r}\sum_{|n|=0}^MC(n)A_0A_1^{(n_1)}
\cdots A_m^{(n_m)}
R_s(\lambda)^{m+1+|n|}d\lambda\nno
&+&\frac{1}{2\pi i}\int_l\lambda^{-p/2-r}Pd\lambda,\eean
where the remainder $P$ is of order $-M-3-2m+|k|$. We will ignore the 
remainder for a moment, and observe that the integrals over $\lambda$ in all 
the other terms may be performed using Cauchy's formula
\ben f^{(b)}(z)=\frac{b!}{2\pi i}\int_C\frac{f(\lambda)}{(\lambda-z)^{b+1}}
d\lambda\een
and the derivative formula
\ben 
\frac{d^b}{d\lambda^b}\lambda^{-p/2-r}=(-1)^b\frac{\Gamma(p/2+b+r)}
{\Gamma(p/2+r)}
\lambda^{-p/2-r-b}.\een
There are two subtle points here. One, is pulling the unbounded operator
$A_0A_1^{(n_1)}\cdots A_m^{(n_m)}$ out of the integral, and the second is
the application of Cauchy's Theorem in the operator setting. The second point 
is handled as in the previous application of Cauchy's Theorem in the 
introduction to section \ref{higgystricks}. The first difficulty is overcome by 
noting that
$$\frac{1}{2\pi i}\int_l \lambda^{-p/2-r}
R_s(\lambda)^{m+1+|n|}d\lambda =(-1)^{m+|n|}
\frac{\Gamma(m+|n|+p/2+r)}{\Gamma(p/2+r)(m+|n|)!}(1+\tD^2+s^2)^{-p/2-r-m-|n|}$$
and that both the integrand and the RHS map our core $\HH_\infty$ into 
itself, so that evaluating our integrals on vectors in $\HH_\infty$, we 
see that we can "push" $A_0A_1^{(n_1)}\cdots A_m^{(n_m)}$ (which is defined 
on all of $\HH_\infty$) through the integral
sign. This yields (modulo the remainder term)
\be B(s)=\sum_{|n|=0}^M(-1)^{m+|n|}\frac{\Gamma(m+|n|+p/2+r)}
{\Gamma(p/2+r)(m+|n|)!}C(n)A_0A_1^{(n_1)}\cdots A_m^{(n_m)}
(1+\tD^2+s^2)^{-p/2-r-m-|n|}.\label{n}\ee
Now, there is $B\in OP^0$ with $A_0A_1^{(n_1)}\cdots A_m^{(n_m)}=
B(1+\tD^2)^{-|n|/2}$ so that, 
\ben A_0A_1^{(n_1)}\cdots A_m^{(n_m)}(1+\tD^2+s^2)^{-|k|/2-|n|/2}
=B(1+\tD^2)^{-|n|/2}(1+\tD^2+s^2)^{-|k|/2-|n|/2}\een
is uniformly bounded independent of $s$. So (again modulo the remainder term)
\be \n B(s)\n_1\leq \sum_{|n|=0}^MC'_n\n
(1+\tD^2+s^2)^{-p/2-r-m-|n|/2+|k|/2}\n_1,\label{label}\ee
and this is finite when $Re(r)+m+|n|/2-|k|/2>0$. For the worst case 
($|n|=0$) we 
obtain the condition in the statement of the lemma.

For the remainder term we have, using Lemmas \ref{uni} and \ref{lambda}, 
\bea \n\int_l\lambda^{-p/2-r}Pd\lambda\n_1&\leq& C_r\int_{-\infty}^\infty
\sqrt{a^2+v^2}^{-p/2-Re(r)}\n R_s(\lambda)^{M/2+3/2+m-|k|/2}\n_1dv\nno
\leq  C_r&\displaystyle{\int_{-\infty}^\infty}&\sqrt{a^2+v^2}^{-p/2-Re(r)}
\sqrt{(1/2+s^2-a)^2+v^2}^{-M/2-3/2-m+|k|/2+(p+\epsilon)/2}dv.\label{r}\eea
Thus the remainder term is trace class for $M>\epsilon+|k|-1-2Re(r)-2m$, 
which can always be arranged.

Finally, adding Equations (\ref{label}) and (\ref{r}), multiplying by 
$s^\alpha$, and integrating over $[0,\infty)$ gives (using Lemma \ref{BIG})
\bean &&\int_0^\infty s^\alpha\n B(s)\n_1ds\leq\sum_{|n|=0}^MC_n'\int_0^\infty 
s^\alpha(1/2+s^2)^{\epsilon-Re(r)-m-|n|/2+|k|/2}ds\nno
&+&C_r\int_0^\infty s^\alpha\int_{-\infty}^\infty\sqrt{a^2+v^2}^{-p/2-Re(r)}
\sqrt{(1/2+s^2-a)^2+v^2}^{-M/2-3/2-m+|k|/2+(p+\epsilon)/2}dvds.\eean
Inspection shows that the non-remainder terms are finite for 
\ben \alpha+2\epsilon-2Re(r)-2m-|n|+|k|<-1,\een
which gives the final statement of the lemma upon considering the worst case, 
$|n|=0$. Using Lemma \ref{intest}, the remainder term is finite for 
\ben M>\alpha-2-2m+|k|+p+\epsilon,\een
which may always be arranged.
\end{proof}

\begin{lemma}\label{reso} 
With the notation as set out at the beginning of this 
subsection and with $R_s(\lambda)=(\lambda-1-\tD^2-s^2)^{-1}$,
$\tilde{R}_s(\lambda)=(\lambda-(1+\tD^2+s^2+s\{\tD,q\}))^{-1}$
we have for  $Re(r)>0$  and any positive integer $M>p-1$:
\bean
&&S\tau(q(1+\tD^2+s^2+s\{\tD,q\})^{-p/2-r})\nno
&=&\sum_{m=1, odd}^{M}s^m
S\tau\left(\frac{1}{2\pi i}\int_l\lambda^{-p/2-r}q\left(R_s(\lambda)
\{\tD,q\}\right)^mR_s(\lambda)d\lambda\right)\nno
&+&s^{M+1}S\tau\left(\frac{1}{2\pi 
i}\int_l\lambda^{-p/2-r}q
\left(R_s(\lambda)\{\tD,q\}\right)^{M+1}
\tilde{R}_s(\lambda)d\lambda\right).\eean
\end{lemma}

\begin{proof}
We will use Cauchy's formula 
\ben f(z)=\frac{1}{2\pi i}\int_C\frac{f(\lambda)}{(\lambda-z)}d\lambda\een
and the resolvent expansion (easily proved by induction on $M$)
\ben \tilde{R}_s(\lambda)=\sum_{m=0}^M 
\left(R_s(\lambda)s\{\tD,q\}\right)^m R_s(\lambda)
+\left(R_s(\lambda)s\{\tD,q\}\right)^{M+1} \tilde{R}_s(\lambda),\een
valid for $\lambda\in l$ 
to expand
$q(1+\tD^2+s^2+s\{\tD,q\})^{-p/2-r}.$  
Since $q$ and $\{\tD,q\}$ are in $OP^0$ by {\bf Observation 3} we can employ
the previous lemmas to first use Cauchy's formula to obtain
\ben
S\tau(q(1+\tD^2+s^2+s\{\tD,q\})^{-p/2-r})=
 S\tau\left(\frac{1}{2\pi i}\int_l
\lambda^{-p/2-r}q(\lambda-(1+\tD^2+s\{\tD,q\}+s^2))^{-1}d\lambda\right)\een
and then apply the resolvent expansion to arrive at
\bean && S\tau\left(q(1+\tD^2+s^2+s\{\tD,q\})^{-p/2-r}\right)\nno
&=&S\tau\left(\frac{1}{2\pi 
i}\int_l\lambda^{-p/2-r}\sum_{m=1,odd}^{M}s^mq\left(R_s(\lambda)
\{\tD,q\}\right)^mR_s(\lambda)d\lambda\right)\nno
&+&s^{M+1}S\tau\left(\frac{1}{2\pi 
i}\int_l\lambda^{-p/2-r}q
\left(R_s(\lambda)\{\tD,q\}\right)^{M+1}
\tilde{R}_s(\lambda)d\lambda\right).\eean
Separating out the remainder term is valid, because $M+1>p>p/2$ 
ensures, by Lemma \ref{yucko}, that it is trace-class. 

We have retained only odd terms in the resolvent expansion, because the even 
terms 
vanish under the supertrace. For if we consider a single term in the sum with 
$k$ 
$\{\tD,q\}$'s, we find (since $\rho^2=1$ and $\rho$  commutes with $\tD$ and 
$\Gamma$, but anticommutes with $q$)
\bean 
&&S\tau\left(\frac{1}{2\pi 
i}\int_l\lambda^{-p/2-r}\rho^2q\left(R_s(\lambda)\{\tD,q\}\right)^k
R_s(\lambda)d\lambda\right)\nno
&&=(-1)^{k+1}S\tau\left(\frac{1}{2\pi 
i}\int_l\lambda^{-p/2-r}\rho 
q\left(R_s(\lambda)\{\tD,q\}\right)^kR_s(\lambda)\rho d\lambda\right)\nno
&&=(-1)^{k+1}S\tau\left(\frac{1}{2\pi 
i}\int_l\lambda^{-p/2-r} 
q\left(R_s(\lambda)\{\tD,q\}\right)^kR_s(\lambda) d\lambda\right),\eean
So if $k$ is even we get zero. This argument does not apply to the remainder 
term 
\ben s^{M+1}S\tau\left(\frac{1}{2\pi 
i}\int_l\lambda^{-p/2-r}q
\left(R_s(\lambda)\{\tD,q\}\right)^{M+1}
\tilde{R}_s(\lambda)d\lambda\right).\een
as $\rho$ neither commutes nor anticommutes with $\tilde{R}_s(\lambda)$.

Now, we may apply the conclusion of Lemma \ref{spaz} to the non-remainder terms 
to see that each integral in the sum is trace-class. Hence we may move the 
trace through the sum and obtain the final statement of the Lemma. We note that 
we cannot push the trace through the integrals in the non-remainder terms as
these integrands are not trace-class: we did {\bf not} apply the actual
pseudodifferential expansion in the proof of Lemma \ref{spaz} to each of these 
terms to obtain trace-class integrands.
\end{proof}
The main result of this subsection is the following Lemma.

\begin{lemma}\label{resolvent} 
Let $N=[p/2]+1$ be the least positive integer strictly greater than p/2. Then
there is a $\delta'$, $0<\delta'<1$
such that  
\bean  sf(\D,u^*\D u){C_{p/2+r}}&=&
\int_0^\infty S\tau(q(1+\tD^2+s^2+s\{\tD,q\})^{-p/2-r})ds\nno
&=& \frac{1}{2\pi i}\sum_{m=1,odd}^{2N-1}\int_0^\infty s^m S\tau\left(\int_l 
\lambda^{-p/2-r}q\left(R_s(\lambda)\{\tD,q\}\right)^mR_s(\lambda)d\lambda\right)
ds\nno
&+&holo,
\eean
where $holo$ is a function of $r$ holomorphic for $Re(r)>-p/2+\delta'/2$.
\end{lemma}

\begin{proof}
Writing $A=\{\tD,q\}$, we have by H\"{o}lder's inequality
\ben \n(R_s(\lambda)A)^{2N}\n_1\leq\left(\n R_s(\lambda)A\n_{2N}\right)^{2N}
\leq \n A\n^{2N}\n R_s(\lambda)\n_{2N}^{2N}=\n A\n^{2N}\n R_s(\lambda)^{2N}\n_1.
\een
Consequently, we can estimate the super-trace of the remainder term
from Lemma \ref{reso} using Lemmas \ref{lambda} and 
\ref{intest}:
\bean 
\int_0^\infty s^{2N}\n\int_l\lambda^{-p/2-r}q\left(R_s(\lambda)\{\tD,q\}
\right)^{2N}\tilde{R}_s(\lambda)d\lambda\n_1 ds\nno
\leq 
C\int_0^\infty s^{2N}\int_{-\infty}^\infty\sqrt{a^2+v^2}^{-p/2-Re(r)}\n 
R_s(\lambda)^{2N}\n_1\n \tilde{R}_s(\lambda)\n dvds\eean
Then  Lemmas \ref{lambda} and \ref{normlambda} give the bound
$$\leq C'\int_0^\infty \int_{-\infty}^\infty s^{2N}
\frac{(a^2+v^2)^{-p/4-Re(r)/2}((1/2+s^2-a)^2+v^2)^{(p+\epsilon)/4-N}}
{((1+s^2-a-sc)^2+v^2)^{1/2}}dvds,$$
where $c=\n\{\tD,q\}\n$.
We now apply Lemma \ref{intest} with $J=2N$, $M=p/2+r$, $A=1$ and 
$K=2N-(p+\epsilon)/2$. A simple 
check shows that the preceding integral is finite for $Re(r)>-N+(1+\epsilon)/2.$
Letting $N = p/2+\delta/2$ where $1\geq \delta >0$ by Observation 1 of this 
Section, we see that the integral is finite for 
$Re(r)>-p/2+(1-\delta +\epsilon)/2$.
Choose $\epsilon<\delta$ and then $\delta^{\prime}=(1-\delta +\epsilon).$ 

To see that the remainder term:
$$F(r)=\frac{1}{2\pi i}\int_0^{\infty}s^{2N}\int_l\lambda^{-p/2-r}S\tau
\left(q(R_s(\lambda)\{\tD,q\})^{2N}\tilde R_s(\lambda)\right)d\lambda ds$$
is defined and holomorphic for $r$ in the half-plane, 
$Re(r)>-p/2+\delta^{\prime}/2,$ we first fix $r.$ Since $\lambda\in l$,
$|\lambda^{-p/2-r}|\leq C_r|\lambda|^{-p/2-Re(r)},$
so the integral converges by Lemma \ref{intest}. We will show
that $F^{\prime}(r)$ is given by the above integral with $\lambda^{-p/2-r}$ 
replaced by $-\lambda^{-p/2-r}Log(\lambda).$ 
To see that this new integral converges, we observe that since $r$ is fixed, 
there exists a fixed $\eta>0$ with $Re(r)-\eta>-p/2+\delta^{\prime}/2,$
and so:
$$|\lambda^{-p/2-r}Log(\lambda)|\leq C |\lambda|^{-p/2-Re(r)+\eta}
|\lambda|^{-\eta}|Log(\lambda)| \leq C^{\prime}|\lambda|^{-p/2-(Re(r)+\eta)}.$$
So, that integral converges by Lemma \ref{intest}. Now, the 
difference quotient $1/h(F(r+h)-F(r))$ is the above integral, with
$\lambda^{-p/2-r}$ replaced by $\lambda^{-p/2-r}(1/h)(\lambda^h -1)$.
We use Taylor's theorem applied to $f(h)=\lambda^{-h}$ about $0$ 
($\lambda$ fixed) to estimate the difference 
$(1/h)(\lambda^h -1)-(-Log\lambda)$.
So, let $C$ be the circle $z=\eta e^{i\theta}$ for 
$0\leq\theta\leq 2\pi$, and consider $h<\eta/2$. By (\cite{A} pp. 125-6):
$$\frac{\lambda^h -1}{h}-(-Log\lambda)=f_2(h)h\;\;
where\;\;f_2(h)=\frac{1}{2\pi i}\int_C\frac{\lambda^{-z}dz}{z^2(z-h)}.$$
Moreover,
$$|f_2(h)|\leq \frac{max\{|\lambda^{-z}|:z\in C\}}{\eta\eta/2}\leq
\frac{2|\lambda|^{\eta} e^{\eta\pi/2}}{\eta^2}=(const)|\lambda|^{\eta}.$$
Finally,
$$\left |\lambda^{-p/2-r}\frac{\lambda^{-\eta}-1}{h}-\lambda^{-p/2-r}
(-Log\lambda)\right |\leq |\lambda^{-p/2-r}|\cdot|f_2(h)h|\leq
C_2|h||\lambda|^{-p/2-Re(r)}.$$
Hence the quotient $1/h(F(r+h)-F(r))$ differs from the formal derivative by
$|h|$ times an absolutely convergent integral with nonegative integrand. 
\end{proof}

\subsection{The Resolvent Cocycle}\label{rescocycle}

At this point it is interesting to perform the `super' bit of the trace, 
so that we 
have an expression which only depends on our original spectral triple 
$(\A,\HH,\D)$. 
The computation begins by recalling the definition of $q$ and $\{\tD,q\}$. Then,
\allowdisplaybreaks\bean &&q\left(R_s(\lambda)\{\tD,q\}\right)^mR_s(\lambda)=
\\&=&i(-1)^{(m-1)/2}\s_3\s_1^m 
\otimes\bma u^*R[\D,u]R[\D,u^*]\cdots 
[\D,u]R & 0\\ 0 & uR[\D,u^*]R[\D,u]\cdots 
[\D,u^*]R\ema.\eean
On the right hand side, by an abuse of notation, we have written 
$R=(\lambda-(1+\D^2+s^2))^{-1}$.
The grading operator is $\Gamma=\s_2\otimes\sigma_3\otimes 1$. 
Since $\s_2\s_3\s_1^m=i1_2$ for $m$ odd, we have (writing $Tr_4$ for the  
operator-valued trace which maps 
$\tilde{\mathcal N}=M_2\otimes M_2\otimes \mathcal N\to \mathcal N$)
\bean &&Tr_4(\Gamma qR_s(\lambda)\{\tD,q\}R_s(\lambda)\cdots\{\tD,q\}
R_s(\lambda))\nno
&=&2(-1)^{(m+1)/2}(u^*R[\D,u]R[\D,u^*]
\cdots [\D,u]R-uR[\D,u^*]R[\D,u]
\cdots [\D,u^*]R).\eean
Consequently, there is a $\delta'$ with $0<\delta'<1$ such that for $r>0$ 
\allowdisplaybreaks\bean &&\int_0^\infty 
S\tau(q(1+\tD^2+s^2+s\{\tD,q\})^{-p/2-r}ds\nno
&=& \frac{1}{2\pi i}\sum_{m=1,odd}^{2N-1}\int_0^\infty s^m 
S\tau\left(q\int_l \lambda^{-p/2-r}
\left(R_s(\lambda)\{\tD,q\}\right)^mR_s(\lambda)d\lambda\right)ds+holo\nno
&=&\frac{1}{2\pi i}\sum_{m=1,odd}^{2N-1}\int_0^\infty s^m 
\frac{1}{2}\tau\left(Tr_4\Gamma q
\left(\int_l \lambda^{-p/2-r}\left(R_s(\lambda)\{\tD,q\}\right)^m
R_s(\lambda)d\lambda\right)\right)ds+holo\nno
&=&\frac{1}{2\pi i}\sum_{m=1,odd}^{2N-1}\int_0^\infty s^m \frac{1}{2}\tau
\left(\int_l \lambda^{-p/2-r}Tr_4\left(\Gamma q\left(R_s(\lambda)\{\tD,q\}
\right)^m
R_s(\lambda)d\lambda\right)\right)ds+holo\nno
&=&\frac{1}{2\pi i}\sum_{m=1,odd}^{2N-1}(-1)^{(m+1)/2}\int_0^\infty s^m\times 
\nno
&\times&\tau\left(\int_l \lambda^{-p/2-r}
\left(u^*R[\D,u]R[\D,u^*]\cdots [\D,u]R-uR
[\D,u^*]R[\D,u]\cdots [\D,u^*]R\right)d\lambda\right)ds\nno
&+&holo\eean
where $holo$ is a function of $r$ holomorphic for $Re(r)>-p/2+\delta'/2$.

$\spadesuit$ This last expression suggests how we might define
a one parameter family of functionals which form a cyclic cocycle 
(up to functionals holomorphic for $Re(r)>-p/2+\delta'/2$)
which we term the {\bf resolvent cocycle}. It will eventuate that 
the Connes-Moscovici residue cocycle may be derived from this resolvent cocycle.
We have been using the notation $R_s(\lambda)$
for the $\tD$ resolvent and $R$ for the $\D$ resolvent. For the remainder
of this section, $R_s(\lambda)=(\lambda-(1+\D^2+s^2))^{-1}$, the $\D$ 
resolvent. $\spadesuit$

\begin{defn}\label{expectation} 
For $m\geq 0$, operators $A_0,...,A_m$, $A_i\in OP^{k_i}$, and 
$2Re(r)>k_0+\cdots+k_m-2m$  define
\ben \la A_0,...,A_m\ra_{m,s,r}=\tau\left(\frac{1}{2\pi i}\int_l 
\lambda^{-p/2-r}
A_0R_s(\lambda)A_1\cdots A_mR_s(\lambda)d\lambda\right).\een
\end{defn}

`Expectations' like these  were first considered by Higson in \cite{hig}. The 
conditions on the orders and on $r$ are sufficient for the trace to be 
well-defined, by Lemma \ref{spaz}. 

\begin{lemma}\label{strick} For any integers $m\geq 0, k\geq 1$
and 
operators $A_0,...,A_m$ with $A_j\in OP^{k_j}$, and $2Re(r)>k+\sum k_j-2m$, 
we may choose $r$ with $Re(r)$ sufficiently large such that
\ben k\int_0^\infty s^{k-1}\la A_0,...,A_m\ra_{m,s,r}ds=
-2\sum_{j=0}^m\int_0^\infty s^{k+1}\la 
A_0,...,A_j, 1, A_{j+1},...,A_m\ra_{m+1,s,r}ds.\een
\end{lemma}
\begin{proof}
For $Re(r)$ sufficiently large, determined by Lemma \ref{spaz}, we have
\bean &&\frac{d}{ds} s^k\frac{1}{2\pi i}\tau\left( \int_l\lambda^{-p/2-r} 
A_0R_s(\lambda)A_1
\cdots A_mR_s(\lambda)d\lambda\right)\nno
&=& ks^{k-1}\frac{1}{2\pi i}\tau\left( \int_l\lambda^{-p/2-r} A_0
R_s(\lambda)A_1
\cdots A_mR_s(\lambda)d\lambda\right)+\nno
&+&\sum_{j=0}^m2s^{k+1}\frac{1}{2\pi i}\tau\left( \int_l\lambda^{-p/2-r} A_0
R_s(\lambda)A_1
\cdots A_jR_s(\lambda)^2A_{j+1}\cdots A_mR_s(\lambda)d\lambda\right).\eean
The fundamental theorem of calculus completes the argument.
\end{proof}

\begin{lemma}\label{cyclicity} For $m\geq 0$, $Re(r)$
sufficiently large, 
$k\geq 1$ and  $A_j\in OP^{(k_j)}$, $j=0,...,m$, we have
\ben \int_0^\infty s^k\la A_0,...,A_m\ra_{m,s,r}ds
=\int_0^\infty s^k\la A_m,A_0,...,A_{m-1}
\ra_{m,s,r}ds.\een
\end{lemma}
\begin{proof} The size of $Re(r)$ is determined by Definition \ref{expectation} 
via 
Lemma \ref{spaz}. We can repeat Lemma \ref{strick} until the integrand of 
\ben\la A_0,...1,...,1,...,A_m\ra_{m+M,s,r}\een
 is trace class. The cyclicity of the 
trace allows us to conclude.
\end{proof}

\begin{lemma}\label{anticomm} For operators $A_0,...,A_m$, $A_j\in OP^{k_j}$, 
$k_j\geq 0$, and $Re(r)$ sufficiently large we have
\bean &&-\la A_0,...,[\D^2,A_j],...,A_m\ra_{m,s,r}\nno
&=&\la A_0,...,A_{j-1}A_j,...,A_m\ra_{m-1,s,r}-\la A_0,...,A_jA_{j+1},...,A_m
\ra_{m-1,s,r},\eean
and for $k\geq 1$
\ben \int_0^\infty s^k\la\D A_0,A_1,...,A_m\ra_{m,s,r}ds=
\int_0^\infty s^k\la A_0,A_1,...,A_m\D\ra_{m,s,r}ds.\een
\end{lemma}
\begin{proof} The first identity follows from observing that 
\ben -[\D^2,A_j]=R_s(\lambda)^{-1}A_j-A_jR_s(\lambda)^{-1},\een
and cancelling neighbouring $R_s(\lambda)$'s. The second follows by applying
Lemma \ref{cyclicity}, then commuting $\D$ past the (hidden) $R_s(\lambda)$
and applying Lemma \ref{cyclicity} again:
\bean\int_0^\infty s^k\la\D A_0,A_1,...,A_m\ra_{m,s,r}ds&=&
\int_0^\infty s^k\la A_m,\D A_0,A_1,...,A_{m-1}\ra_{m,s,r}ds\nno
=\int_0^\infty s^k\la A_m\D, A_0,A_1,...,A_{m-1}\ra_{m,s,r}ds&=&
\int_0^\infty s^k\la A_0,A_1,...,A_m\D\ra_{m,s,r}ds\eean
\end{proof}
Suspecting that the spectral flow is given by pairing a cocycle with the 
Chern character of a unitary, we remove the normalisation coming from 
$Ch_m(u)$ from our resolvent formula to define a cocycle. The factor of 
$\sqrt{2\pi i}$ is for compatability with the Kasparov product, \cite{Co4}.
\begin{defn} Let $\mathcal C(m)$ denote the constant
$\frac{-2\sqrt{2\pi i}}{\Gamma((m+1)/2)}.$
Then, for $Re(r)>-m/2+1/2$ and $da=[\D,a]$ we define 
$\phi_m^r:\A^{m+1}\to{\C}$ by 
\ben \phi_m^r(a_0,...,a_m)=\mathcal C(m)
\int_0^\infty s^m\la a_0,da_1,..,da_m\ra_{m,s,r}ds.\een
\end{defn}
 
By Lemma \ref{spaz} the condition on $r$ ensures that the integral converges.
We note that this constant $\mathcal C(m)$ is distinct from $C(k)$ which takes
a multi-index $k$ as its argument. The next result captures the main new 
idea of this Section.

\begin{prop}\label{cocycle}  
For $p\geq 1$ the collection of functionals $\phi^r=\{\phi_m^r\}_{m=1}^{2N-1}$, 
$m$ odd, is such that 
\be (B\phi^r_{m+2}+b\phi^r_m)(a_0,...,a_{m+1})=0\ \ \ m=1,3,...,2N-3,\ \ \ 
(B\phi^r_1)(a_0)=0\label{coc}\ee
where the $a_i\in\A$. Moreover, there is a $\delta'$, $0<\delta'<1$
such that  
$b\phi_{2N-1}^r(a_0,...,a_{2N})$ is a holomorphic function of $r$ for 
$Re(r)>-p/2+\delta'/2$.
\end{prop}
\begin{proof} 
We start with the computation of the coboundaries of the $\phi_m^r$. The
definition 
of the operator $B$ and $\phi^r_{m+2}$ gives
\bean (B\phi^r_{m+2})(a_0,...,a_{m+1})&=&\sum_{j=0}^{m+1}
\phi^r_{m+2}(1,a_j,...,a_{m+1},a_0,...,a_{j-1})\nno
&=&\sum_{j=0}^{m+1}\mathcal C(m+2)\int_0^\infty s^{m+2}
\la 1,[\D,a_j],...,[\D,a_{j-1}]\ra_{m+2,s,r}ds.\eean
Using Lemma \ref{cyclicity} and Lemma \ref{strick}, this is equal to
\bean &=&\sum_{j=0}^{m+1}\mathcal C(m+2)\int_0^\infty s^{m+2}
\la [\D,a_0],...,[\D,a_{j-1}],1,[\D,a_j],...,[\D,a_{m+1}]\ra_{m+2,s,r}ds\nno
&=&-\mathcal C(m+2)\frac{(m+1)}{2}
\int_0^\infty s^{m}\la [\D,a_0],...,[\D,a_{m+1}]\ra_{m+1,s,r}ds,\eean

We observe at this point that $\mathcal C(m+2)(m+1)/2=\mathcal C(m)$, 
using the functional equation for the Gamma function.

Next we write $[\D,a_0]=\D a_0-a_0\D$ and anticommute the second $\D$ through 
the remaining 
$[\D,a_j]$ using $\D[\D,a_j]+[\D,a_j]\D=[\D^2,a_j]$. This gives us 
\bea &&(B\phi^r_{m+2})(a_0,...,a_{m+1})\nno
&=& -\mathcal C(m)\int_0^\infty s^{m}
\la \D a_0-a_0\D,[\D,a_1],...,[\D,a_{m+1}]\ra_{m+1,s,r}ds\nno
&=&-\mathcal C(m)\int_0^\infty s^{m}
\left(\la \D a_0,[\D,a_1],...,[\D,a_{m+1}]\ra_{m+1,s,r}\right.\nno
&-&\left.\la  a_0,[\D,a_1],...,[\D,a_{m+1}]\D\ra_{m+1,s,r}\right)ds\nno
&-&\mathcal C(m)\int_0^\infty s^m \sum_{j=1}^{m+1}(-1)^{j}
\la a_0,[\D,a_1],...,[\D^2,a_j],...,[\D,a_{m+1}]\ra_{m+1,s,r}ds,\nno
&=&-\mathcal C(m)\int_0^\infty s^m \sum_{j=1}^{m+1}(-1)^{j}
\la a_0,[\D,a_1],...,[\D^2,a_j],...,[\D,a_{m+1}]\ra_{m+1,s,r}ds,
\label{bigbee}\eea
where the last line follows from the second identity of Lemma \ref{anticomm}.
Observe that for $\phi_1^r$ we have
$$ (B\phi_1^r)(a_0)=\frac{\mathcal C(1)}{2\pi i}\int_0^\infty s
\tau\left( \int_l\lambda^{-p/2-r}R_s(\lambda)[\D,a_0]R_s(\lambda)
d\lambda\right)ds=0,$$
by an easy variant of the argument in Lemma \ref{anticomm}.
We now compute the Hochschild coboundary of $\phi_m^r$. From the 
definitions we have
\bean (b\phi_m^r)(a_0,...,a_{m+1})&=&\phi_m^r(a_0a_1,a_2,...,a_{m+1})\nno
&+&\sum_{i=1}^m(-1)^i\phi^r_m(a_0,...,a_ia_{i+1},...,a_{m+1})\nno
&+&\phi^r_m(a_{m+1}a_0,a_1,...,a_m)\nno
&=& \mathcal C(m)\int_0^\infty s^m\left(\la a_0a_1,[\D,a_2],...,[\D,a_{m+1}]
\ra_{m,s,r}+\right.\nno
&+&\left.\sum_{i=1}^m(-1)^i\la a_0,[\D,a_1],...,a_i[\D,a_{i+1}]+[\D,a_i]a_{i+1},
...,[\D,a_{m+1}]
\ra_{m,s,r}+\right.\nno
&+&\left.\la a_{m+1}a_0,[\D,a_1],...,[\D,a_m]\ra_{m,s,r}\right)ds.\eean
We now reorganise the terms so that we can employ the first identity of Lemma 
\ref{anticomm}. So
\bea\label{littlebee} &&(b\phi_m^r)(a_0,...,a_{m+1})\nno
&=&\mathcal C(m)\int_0^\infty s^m
\left(\left(\la a_0a_1,[\D,a_2],...,[\D,a_{m+1}]\ra_{m,r,s}-\la a_0,a_1[\D,a_2],
...,[\D,a_{m+1}]\ra_{m,r,s}
\right)-\right.\nno
&-& \mathcal C(m)\int_0^\infty s^m
\left(\la a_0,[\D,a_1]a_2,...,[\D,a_{m+1}]\ra_{m,r,s}-\la a_0,[\D,a_1],a_2
[\D,a_3],...,[\D,a_{m+1}]\ra_{m,r,s}
\right)+\nno
&&\vdots\nno
&-&\mathcal C(m)\int_0^\infty s^m
\left.\left(\la a_0,[\D,a_1],...,[\D,a_{m}]a_{m+1}\ra_{m,r,s}-
\la a_{m+1}a_0,[\D,a_1],...,[\D,a_{m}]\ra_{m,r,s}
\right)\right)ds\nno
&=&\sum_{j=1}^{m+1}(-1)^{j}
\mathcal C(m)\int_0^\infty s^m\la a_0,[\D,a_1],...,[\D^2,a_j],...,[\D,a_{m+1}]
\ra_{m+1,r,s}ds\nno
\eea
where to get the term with $[\D^2,a_{m+1}]$ we have used Lemma \ref{cyclicity}. 

For $m=1,3,5,...,2N-3$ comparing Equations (\ref{littlebee}) and (\ref{bigbee}) 
now shows that 
\ben (B\phi_{m+2}^r+b\phi_m^r)(a_0,...,a_{m+1})=0.\een
So we just need to check the claim that $b\phi_{2N-1}^r$  is holomorphic for 
$Re(r)>-p/2+\delta'$ for some suitable $\delta'$. 
 By  Lemma \ref{spaz}
\be \int_0^\infty s^m\la a_0,[\D,a_1],...,[\D^2,a_j],...,[\D,a_{m+1}]
\ra_{m+1,r,s}ds\label{ted}\ee
is finite when
\ben Re(r)>\frac{1+m+2\epsilon+1}{2}-(m+1)=
\frac{1-p}{2}+\frac{2\epsilon-1-m+p}{2}.\een
For $m=2N-1=2[p/2]+1$ this reduces to
\ben Re(r)>\frac{1-p}{2}+\frac{2\epsilon-2[p/2]-2+p}{2}.\een
As $p/2-1<[p/2]\leq p/2$, we see that $0>-2[p/2]-2+p\geq-2$, and we may always 
find an $\epsilon>0$ so that $2\epsilon-2[p/2]-2+p<0$. The proof that 
(\ref{littlebee}) is holomorphic for $m=2N-1$ is similar to the 
analyticity proof in 
Lemma \ref{resolvent}.
\end{proof}

Thus we have an odd $(b,B)$ cochain $(\phi^r_m)_{m=1,odd}^{2N-1}$ with values 
in the functions 
which are holomorphic in a half-plane. Moreover, modulo those functions  
holomorphic in the half-plane $Re(r)>-p/2-\delta$, our resolvent cochain is a 
cocycle. This, together with Lemma \ref{resolvent}, 
actually proves Part $1)$ of Theorem \ref{SFLIT}.

\section{The residue cocycle}\label{residue}

This Section proves Theorem \ref{SFLIT}. It is organised into 3 subsections.
In the first of these we begin with the resolvent expansion of the 
spectral flow formula
Lemma \ref{resolvent}. We use the pseudodifferential calculus to derive 
from the resolvent expansion a new 
expression stated as Proposition 8.1.
This leaves us with a formula for spectral flow that involves an
integral over the parameter $s$. By integrating out the $s$ dependence in the 
formula of Proposition 8.1
we find in Subsection 8.2 (Proposition 8.2) a spectral flow
formula which involves a sum of zeta functions.
One immediately recognises that individual
terms in this formula may be obtained from 
our resolvent cocycle of Section 7 by using the pseudodifferential
calculus. Thus, from the resolvent cocycle
we derive in the final subsection, in Theorem 8.4, the residue cocycle. 
Our final formula for
the spectral flow follows immediately by evaluating the residue cocycle on 
$Ch_*(u^*)$.

\subsection{Pseudodifferential Expansion of the Spectral Flow}\label{psiexpan}

The aim of this subsection is to establish just one formula
which is summarised in the following result. Recall that $N=[p/2]+1$.

\begin{prop}
There is a $\delta'$, $0<\delta'<1$ such that 
\bean
&&sf(\D,u^*\D u){C_{p/2+r}}=\int_0^\infty 
S\tau(q(1+\tD^2+s^2+s\{\tD,q\})^{-p/2-r})ds\nno
&=& \frac{1}{2\pi i}\sum_{m=1,odd}^{2N-1}
\int_0^\infty \!\!s^m S\tau\!\left(\int_l 
\lambda^{-p/2-r}\sum_{|k|=0}^{2N-1-m}C(k)q\{\tD,q\}^{(k_1)}\cdots
\{\tD,q\}^{(k_m)}R
_s(\lambda)^{m+1+|k|}d\lambda\right)ds\nno
&&\qquad\qquad\qquad\qquad\qquad\qquad\qquad\qquad\qquad\qquad\qquad\qquad
\qquad\qquad\qquad\qquad +holo\nno
&=&\sum_{m=1,odd}^{2N-1}\sum_{|k|=0}^{2N-1-m}(-1)^{|k|+m}C(k)
\frac{\Gamma(p/2+r+m+|k|)}{\Gamma(p/2+r)(|k|+m)!}\times\nno
&\times &\int_0^\infty 
s^mS\tau\left(q\{\tD,q\}^{(k_1)}\cdots\{\tD,q\}^{(k_m)}
(1+\tD^2+s^2)^{-(p/2+r+m+|k|)}\right)ds+holo,\eean
where $holo$ is a function of $r$ holomorphic for $Re(r)>-p/2+\delta'/2.$ 
Consequently the sum of functions on the right hand side has an analytic 
continuation to a deleted neighbourhood of $r=(1-p)/2$ (given by the left 
hand side) with at worst a simple pole at $r=(1-p)/2$.
\end{prop}

\begin{proof} 
The proof starts by applying 
the pseudodifferential expansion to each of the terms of 
the resolvent expansion of Lemma \ref{resolvent}. Thus, modulo the holo
of Lemma \ref{resolvent}:
$$\int_0^\infty S\tau(q(1+\tD^2+s^2+s\{\tD,q\})^{-p/2-r})ds$$
\bean &=& \frac{1}{2\pi i}\sum_{m=1,odd}^{2N-1}
\int_0^\infty s^m S\tau\{\int_l 
\lambda^{-p/2-r}\sum_{|k|=0}^{2N-1-m}C(k)q\{\tD,q\}^{(k_1)}\cdots
\{\tD,q\}^{(k_m)}R
_s(\lambda)^{m+1+|k|}d\lambda\}ds\nno
&+& \frac{1}{2\pi i}\sum_{m=1,odd}^{2N-1}
\int_0^\infty s^m S\tau\left(\int_l 
\lambda^{-p/2-r}P_{m}(s,\lambda)d\lambda\right)ds.\eean
Here we are considering multi-indices $k$ with $m$ terms, with $|k|=k_1+\cdots 
+k_m$, and using Lemma \ref{uni} with $M=2N-1-m$ to give us that
\ben P_{m}(s,\lambda)(\lambda-(1+\tD^2+s^2))^{N+(m+2)/2}\een
is bounded where the bound is uniform in $s$ and $\lambda$. 
We refer to the last expression in this formula as
the error term. We want to compute the trace norm of the error term, and 
show that it is integrable in $s$. We have by Lemma \ref{lambda}
\bean &&\n \int_0^\infty s^m\int_l\lambda^{-p/2-r}
P_{m}d\lambda ds\n_1\nno
&\leq&C\int_0^\infty s^m\int_{-\infty}^\infty\sqrt{a^2+v^2}^{-p/2-r}\n
(\lambda-(1+\tD^2+s^2))^{-N-(m+2)/2}\n_1dvds\nno
&\leq&C'\int_0^\infty s^m\int_{-\infty}^\infty\sqrt{a^2+v^2}^{-p/2-r}
\sqrt{(1/2+s^2-a)^2+v^2}^{-N-(m+2)/2+(p+\epsilon)/2}dvds.\eean
We now apply Lemma \ref{intest}. Since $N\geq p/2+\delta/2$ and $0<\epsilon<1$,
this tells us that the remainder term is finite for $r>-p/2+(1-\delta)/2$.
That the error term is holomorphic as a function of $r$ is proved exactly
as in the proof of analyticity in Lemma \ref{resolvent}.
We now argue as in the proof of Lemma \ref{spaz} using the Cauchy integral 
formula
\ben f^{(n)}(z)=\frac{n!}{2\pi 
i}\int_C\frac{f(\lambda)}{(\lambda-z)^{n+1}}d\lambda\een
and the derivative formula\ben 
\frac{d^b}{d\lambda^b}\lambda^{-p/2-r}=(-1)^b
\frac{\Gamma(p/2+b+r)}{\Gamma(p/2+r)}\lambda^{-p/2-r-b}.\een
Applying the Cauchy formula produces the expression in the statement of the 
Proposition.
\end{proof} 
We now have an expression which is a finite sum of terms, 
and which can be used to compute the residues at $r=(1-p)/2$.

\subsection{Eliminating the $s$ dependence in Proposition 8.1}\label{sint}

In this subsection we  integrate out the $s$ dependence
in the formula of Proposition 8.1 to obtain the following new expression.
\begin{prop}\label{killextraterms}
There is a $\delta'$, $0<\delta'< 1$ such that:
\bean &&sf(\D,u^*\D u)C_{p/2+r}\nno
 &=&\sum_{m=1,odd}^{2N-1}\sum_{|k|=0}^{2N-1-m}
C(k)(-1)^{|k|+m}\frac{\Gamma((m+1)/2)\Gamma(p/2+r+|k|+(m-1)/2)}
{2(|k|+m)!\Gamma(p/2+r)}\times\nno
&\times&S\tau\left(q\{\tD,q\}^{(k_1)}\cdots\{\tD,q\}^{(k_m)}
(1+\tD^2)^{-(p/2+r+|k|+(m-1)/2)}\right)+holo,\eean
where $holo$ is a function of $r$ holomorphic for $Re(r)>-p/2+\delta'/2$. 
Consequently the sum of functions on the right hand side has an analytic 
continuation to a deleted neighbourhood of $r=(1-p)/2$ (given by the left 
hand side) with at worst a simple pole at $r=(1-p)/2$.
Moreover, if $[p]=2n$ is even, each of the top terms with $|k|=2N-1-m$
are holomorphic at $r=(1-p)/2$, including the one term with $m=2N-1$.
\end{prop}

\noindent{\bf Remark} 
We started with a spectral flow formula, Equation (\ref{uberformula}), given 
by an integral defined only for $r>0$ (holomorphic for $Re(r)>1/2$ by the 
comments following Equation (\ref{uberformula})).
The above 
result is telling us that the
 sum of zeta functions differs from the left hand side by a 
function which is holomorphic for 
$Re(r)>-p/2+\delta'/2$. Thus this sum of zeta functions has only a 
single simple pole in this same region $Re(r)>-p/2+\delta'$ and the residue at
$r=(1-p)/2$ is the spectral flow with no need to invoke
the isolated spectral dimension hypothesis (although of course we cannot 
conclude that the individual terms are meromorphic in this region). This 
proves part $2)$ of Theorem \ref{SFLIT}.

\begin{proof}
We first need to interchange the $s$ integral in Proposition 8.1 with the 
supertrace. To this end we show that for $Re(r)>-m-|k|/2$  
the function $A_{k,m}:{\R}_+\to\LL^1(\tilde{\mathcal N})$ given by
\ben 
A_{k,m}(s)=q\{\tD,q\}^{(k_1)}\cdots\{\tD,q\}^{(k_m)}
(1+\tD^2+s^2)^{-p/2-r-m-|k|}\een
is continuous. Fix $s_0$ and let $B=q\{\tD,q\}^{(k_1)}\cdots\{\tD,q\}^{(k_m)}$.
Since $B\in OP^{|k|}$, by Lemma \ref{OPs} $B=A(1+\tD^2+s_0^2)^{|k|/2}$ where
$A\in OP^0$. Then, for $t\geq 0$ the resolvent equation gives:
\bean &&
\n B(1+\tD^2+s_0^2)^{-p/2-r-m-|k|} - B(1+\tD^2+t^2)^{-p/2-r-m-|k|}\n_1\nno
&\leq& \n A\n\cdot\n (1+\tD^2+s_0^2)^{-p/2-r-m-|k|/2}\n_1\cdot |t-s|\cdot
\n (1+\tD^2+t^2)^{-p/2-r-m-|k|}\n\nno
&\leq&\n A\n\cdot\n (1+\tD^2+s_0^2)^{-p/2-r-m-|k|/2}\n_1\cdot |t-s|\cdot 1.
\eean
So, we can employ the Bochner integral to get: 
$$ \int_0^\infty s^mS\tau(A_{k,m}(s))ds=S\tau(\int_0^\infty s^mA_{k,m}(s)ds).$$
The $s$ integral on the right side can now be performed using the following
Laplace transform computation. This  
requires an interchange of order of integration, which follows from 
a small variant of  the argument of Lemma 9.1 of \cite{CP2}. So
\bean &&\int_0^\infty s^m(1+s^2+\tD^2)^{-(|k|+m+p/2+r)}ds\\
  &=&\frac{1}{\Gamma(p/2+r+|k|+m)}\int_0^\infty s^m\int_0^\infty 
u^{|k|+m+p/2+r-1}e^{-(1+\tD^2)u}e^{-s^2u}duds \\
&=&\frac{1}{\Gamma(|k|+m+p/2+r)}\int_0^\infty u^{|k|+m+p/2+r-1}\left
(\int_0^\infty s^me^{-s^2u}ds\right)e^{-(1+\tD^2)u}du\\
&=&\frac{\Gamma((m+1)/2)}{2\Gamma(|k|+m+p/2+r)}\int_0^\infty 
u^{|k|+(m-1)/2+p/2+r-1}e^{-(1+\tD^2)u}du \\
&=&\frac{\Gamma((m+1)/2)\Gamma(|k|+(m-1)/2+p/2+r)}{2\Gamma(|k|+m+p/2+r)}
(1+\tD^2)^{-(|k|+(m-1)/2+p/2+r)} .\eean
Observe the factor of two introduced by this computation. Since both the 
integrand and the final result map the core subspace 
$\tilde{\mathcal H}_{\infty}$ into itself, we can ``push'' the operator
$q\{\tD,q\}^{(k_1)}\cdots\{\tD,q\}^{(k_m)}$ through the integral sign.

Applying the above 
calculations we now obtain
a formula in which the $s$ and $\lambda$ dependence has 
been integrated out. This yields 
\bean &&\int_0^\infty S\tau\left(q(1+\tD^2+s^2+s\{\tD,q\})^{-p/2-r}\right)ds\nno
&=&\sum_{m=1,odd}^{2N-1}\sum_{|k|=0}^{2N-1-m}C(k)(-1)^{|k|+m}
\frac{\Gamma((m+1)/2)\Gamma(p/2+r+|k|+(m-1)/2)}
{2(|k|+m)!\Gamma(p/2+r)}\times\nno
&\times&S\tau\left(q\{\tD,q\}^{(k_1)}\cdots
\{\tD,q\}^{(k_m)}(1+\tD^2)^{-p/2-r-|k|-(m-1)/2}\right)+holo,\eean
where $holo$ is a function of $r$ holomorphic for 
$Re(r)>-p/2+\delta'/2$.
We now observe that
\bean &&\left|S\tau\left(q\{\tD,q\}^{(k_1)}\cdots\{\tD,q\}^{(k_m)}
(1+\tD^2)^{-p/2-r-|k|-(m-1)/2}\right)\right|\nno
&&\quad\leq C\n(1+\tD^2)^{-p/2-r-|k|/2-(m-1)/2}\n_1.\eean
This is finite when
\ben r>\frac{1}{2}-\frac{m}{2}-\frac{|k|}{2}=\frac{1-p}{2}+\frac{p-m-|k|}{2}.
\een
So whenever $m+|k|>p$, we obtain a term holomorphic at $r=(1-p)/2$ which may be 
discarded. 
Thus for $[p]=2n$ we see $[p/2]=n$ and $2N-1=2([p/2]+1)-1=2[p/2]+1=2n+1>p$. So 
for $m=2N-1$ (and so $|k|=0$), $m+|k|>p$ and this term is holomorphic at 
$r=(1-p)/2$. Similarly when $|k|=2N-1-m$ for any $m=1,3,...,2N-1$, we have 
$m+|k|=2N-1>p$.
\end{proof}

Observe that at this point we have almost proved part $2)$ of 
Theorem \ref{SFLIT}. The only outstanding item is the precise form of the 
constants, and these are identified in the next section.

\subsection{The Residue Cocycle}\label{residuecocycle}

In this subsection we will prove Theorem 4.2.
We have now come to the point where we need to assume that
our semi-finite spectral triple 
$(\A,\HH,\D)$ has isolated spectral dimension.
So this means we may analytically continue our zeta functions to 
a deleted neighbourhood of the critical point. 
We denote the analytic continuation of a function analytic in a 
right half plane to a deleted neighbourhood of the critical point by putting
the function in boldface.
Thus define the functionals for each integer $j\geq 0$:
\bean 
&&S\tau_j(q\{\tD,q\}^{(k_1)}\cdots\{\tD,q\}^{(k_m)}(1+\tD^2)^{-(|k|+m/2)})=
\nno
&=&res_{r=(1-p)/2}(r-(1-p)/2)^j
{\bf S\tau(q\{\tD,q\}^{(k_1)}\cdots\{\tD,q\}^{(k_m)}
(1+\tD^2)^{-(r-(1-p)/2+|k|+m/2)})}.\eean
Analogously we can define $\tau_j$ in terms of the trace $\tau$. Observe that 
taking residues and performing the `super' bit of the trace commute.
Thus we have the Laurent expansion 
\bean 
&&{\bf S\tau(q\{\tD,q\}^{(k_1)}\cdots\{\tD,q\}^{(k_m)}
(1+\tD^2)^{-(r-(1-p)/2+|k|+m/2)})}
=\nno
&=&\sum_{j\geq 0}(r-(1-p)/2)^{-j-1}S\tau_j(q\{\tD,q\}^{(k_1)}\cdots
\{\tD,q\}^{(k_m)}(1+\tD^2)^{-(|k|+m/2)})+holo.\eean
Here $holo$ is a function of $r$ holomorphic for 
$Re(r)>(1-p)/2-\delta$.

Now we start from the formula of Proposition 8.2
and take residues at $r=(1-p)/2$ of both sides.
By our assumption that the zeta functions analytically continue to a deleted
neighbourhood of this critical point (and using the boldface 
notational convention as above)
\bean sf(\D,u^*\D u) 
&=&\sum_{m=1,odd}^{2N-1}\sum_{|k|=0}^{2N-1-m}(-1)^{|k|+m}
\frac{\Gamma((m+1)/2)}{2}\alpha(k)
\times\eean
\vspace{-.04in}
$$res\left(\frac{\Gamma(|k|+(p+m-1)/2+z)}{\Gamma(p/2+z)}
{\bf S\tau(q\{\tD,q\}^{(k_1)}\cdots\{\tD,q\}^{(k_m)}
(1+\tD^2)^{-(z-(1-p)/2+|k|+m/2)})}\right)$$
Now we have,  writing the integer $|k|+(m-1)/2$ as $h=|k|+(m-1)/2$,
\bean \frac{\Gamma(p/2+h+z)}{\Gamma(p/2+z)}&=&
\frac{\Gamma((z-(1-p)/2)+h+1/2)}{\Gamma((z-(1-p)/2)+1/2)}\nno
=\prod_{j=0}^{h-1}\frac{\Gamma((z-(1-p)/2)+j+1+1/2)}{\Gamma((z-(1-p)/2)+j+1/2)}
&=&\prod_{j=0}^{h-1}((z-(1-p)/2)+j+1/2)\nno
&=&\sum_{j=0}^{h}(z-(1-p)/2)^j\s_{h,j}.\eean
Recall that the $\s_{h,j}$'s are the symmetric functions of the 
half-integers $1/2,3/2,...,h-1/2$. 
So finally, we obtain the following formula for $sf(\D,u^*\D u)$
\bean 
&&\sum_{m=1,odd}^{2N-1}\sum_{|k|=0}^{2N-1-m}(-1)^{|k|+m}
\frac{\Gamma((m+1)/2)}{2}\alpha(k)\sum_{j=0}^h\s_{h,j}\times\nno
&\times&res\left((z-(1-p)/2)^j{\bf S\tau
(q\{\tD,q\}^{(k_1)}\cdots\{\tD,q\}^{(k_m)}(1+\tD^2)^{-(z-(1-p)/2+|k|+m/2)})}
\right)\nno
&=&\sum_{m=1,odd}^{2N-1}\sum_{|k|= 0}^{2N-1-m}(-1)^{|k|+m}
\frac{\Gamma((m+1)/2)}{2}\alpha(k)
\times\nno
&\times&\sum_{j=0}^h\s_{h,j}S\tau_j(q\{\tD,q\}^{(k_1)}\cdots
\{\tD,q\}^{(k_m)}(1+\tD^2)^{-(|k|+m/2)})\nno
&=&\sum_{m=1,odd}^{2N-1}\sum_{|k|=0}^{2N-1-m}(-1)^{(m+1)/2+|k|}
\frac{\Gamma((m+1)/2)
}{2}\alpha(k)\times\nno
&\times&\sum_{j=0}^{h}\s_{h,j}\tau_j\left((u[\D,u^*]^{(k_1)}\cdots
[\D,u^*]^{(k_m)}-u^*
[\D,u]^{(k_1)}\cdots[\D,u]^{(k_m)})(1+\D^2)^{-|k|-m/2}\right).\eean
The last line follows from converting the super trace into an ordinary trace. 
The normalisation of $1/2$ for the super trace has been cancelled by a trace
over the $2\times 2$ identity matrix, so we still have a factor of a $1/2$
which arose from the $s$-integral.
 
To understand this formula in terms of cyclic (co)homology and Chern characters,
we show that our spectral flow formula is obtained by pairing a cyclic cocycle 
with the Chern character of a unitary. 
We note that our resolvent cocycle $\phi^r$ pairs with normalised
chains so that by Lemma \ref{unitary},
$$\phi^r(Ch_*(u))=-\phi^r(Ch_*(u^{*}))$$
modulo functions holomorphic for $Re(r)>(1-p)/2-\delta$.

In the next theorem we introduce the functionals
that form the residue cocycle and
from which we obtain part $3)$ of Theorem \ref{SFLIT}.

\begin{thm}\label{grandfinale}
Assume that $(\A,\HH,\D)$ is an odd $QC^\infty$ spectral triple with isolated 
spectral dimension spectrum $p\geq 1$. For $m$ odd, define functionals 
$\phi_m$ by 
\ben \phi_m(a_0,...,a_m)=\sqrt{2\pi 
i}\sum_{|k|=0}^{2N-1-m}\!\!(-1)^{|k|}\alpha(k)\!
\sum_{j=0}^{h}\s_{h,j}\tau_j\left(a_0[\D,a_1]^{(k_1)}\cdots[\D,a_m]
^{(k_m)}(1+\D^2)^{-|k|-m/2}\right),\een
where $h=|k|+(m-1)/2$.
Then $\phi=(\phi_m)$ is a $(b,B)$-cocycle and
\ben  sf(\D,u^*\D u)=\frac{1}{\sqrt{2\pi 
i}}\sum_{m=1}^{2N-1}\langle\phi_m,Ch_m(u)\rangle\een
\end{thm}

\begin{proof}
The proof is implicit in what has gone before. However we
set out the brief direct argument that
the $\phi_m$ define a cocycle as this is by far the most
difficult part of \cite{CM}.

We know that the resolvent cocycle is given by
\ben \phi_m^r(a_0,...,a_m)=
\mathcal C(m)\int_0^\infty s^m\la a_0,[\D,a_1],..,[\D,a_m]\ra_{m,s,r}ds.\een
We now employ the same arguments as in the proof of Propositions 8.1 to obtain 
{\em modulo functions of $r$ holomorphic for $Re(r)>(1-p)/2-\delta$}
\bean &&\phi_m^r(a_0,...,a_m)=\frac{\mathcal C(m)}{2\pi i}\sum_{|k|=0}^{2N-1-m}
C(k)\int_0^\infty s^m\tau(\int_\ell \lambda^{-p/2-r}A_kR_s(\lambda)^{|k|+m+1}
d\lambda)ds,\eean
where $A_k= a_0[\D,a_1]^{(k_1)}...[\D,a_m]^{(k_m)}$ are fixed (for this 
discussion) operators of order $|k|$. We find, by the same 
proof as Proposition 8.2,
\bean &&\phi_m^r(a_0,...,a_m)\nno
&=&\sqrt{2\pi i}\sum_{|k|=0}^{2N-1-m}
\frac{2C(k)(-1)^{|k|}\Gamma(p/2+r+|k|+m)}
{\Gamma((m+1)/2)\Gamma(p/2+r)(|k|+m)!}\int_0^\infty 
s^m\tau(A_k(1+s^2+\D^2)^{-p/2-r-|k|-m})ds\nno
&=&\sqrt{2\pi i}\sum_{|k|=0}^{2N-1-m}
\frac{(-1)^{|k|}\alpha(k)\Gamma(p/2+r+|k|+(m-1)/2)}
{\Gamma(p/2+r)}\tau(A_k(1+\D^2)^{-p/2-r-|k|-(m-1)/2})\nno
&=&\sqrt{2\pi i}\sum_{|k|=0}^{2N-1-m}(-1)^{|k|}\alpha(k)
\sum_{j=0}^h(r-(1-p)/2)^j\sigma_{h,j}
\tau(A_k(1+\D^2)^{-p/2-r-|k|-(m-1)/2}),\eean
where $h=|k|+(m-1)/2$ and $\alpha(k)=\frac{1}{k_1!k_2!\cdots 
k_m!(k_1+1)(k_1+k_2+2)\cdots(|k|+m)}$. 

Taking the residue of the analytic continuations at the critical point 
$$res_{r=(1-p)/2}{\bf \Phi_m^r(a_0,...,a_m)}=\sqrt{2\pi i}
\sum_{|k|=0}^{2N-1-m}\alpha(k)(-1)^{|k|}
\sum_{j=0}^h\sigma_{h,j}\tau_j(A_k(1+\D^2)^{-|k|-m/2}).$$
That is 
$$res_{r=(1-p)/2} {\bf\Phi_m^r(a_0,...,a_m)}=\phi_m(a_0,...,a_m)$$
The proof that $b\phi_{m-2}+ B\phi_m=0$ now follows from the
fact that the $\phi^r_m$ satisfy this linear relation
modulo functions holomorphic for $Re(r)>(1-p)/2-\delta$
so the analytic continuations of the zeta functions
satisfy a similar linear relation and so do the residues, since residues depend 
linearly on the analytic function.
\end{proof}

Part $3)$ of Theorem \ref{SFLIT} follows
from the fact that the $\phi_j$ form a cocycle, and Lemma \ref{unitary}.
We have (with $(y)$ representing the boundary chain in the proof of Lemma 
\ref{unitary})
\bean  sf(\D,u^*\D u)=\frac{1}{2\sqrt{2\pi i}}
\sum_{m=1}^{2N-1}\langle\phi_m,2Ch_m(u)+(b+B)(y)\rangle
=\frac{1}{\sqrt{2\pi i}}\sum_{m=1}^{2N-1}\langle\phi_m,Ch_m(u)\rangle 
\eean
which is just a rewriting of the formula of part $3)$ of Theorem \ref{SFLIT}.

Notice that we have, by this result, 
the correct normalisation for the pairing between cyclic cohomology and 
$K_1({\mathcal A})$.
Observe that since $h=|k|+(m-1)/2$, and $m\leq 2N-1$, $|k|\leq 2N-m-1$, 
no matter 
what order the singularity of the trace is, we need only consider the first $h$ 
terms in the principal part of the Laurent expansion, and
\ben h\leq 2N-m-1+(m-1)/2\leq 2N-(m+3)/2\leq 2N-2\leq 2(p/2+1)-2=p.\een
Observe that this is the renormalised version of Connes-Moscovici's formula. 
This 
is, essentially, because we started from a  scale invariant formula, whereas 
Connes-Moscovici started from the JLO cocycle. The upshot is that they 
obtained a  formula with infinitely many terms, 
and then used the easily determined behaviour of the functionals 
$\tau_j$ under change of scale to obtain counterterms. 
Alternatively, they  replace $\zeta\to 
res_{s=0}\Gamma(h+s+1/2)\zeta(s)$ by the functional 
\ben \zeta\to res_{s=0}\frac{\Gamma(1/2)\Gamma(h+s+1/2)}{\Gamma(s+1/2)}\zeta(s).
\een
thus allowing the removal of the hypothesis of finiteness of 
the dimension spectrum. It is not clear to us whether one can perform this 
manoeuvre in the middle of their proof, or whether it is necessary for them to 
first consider the unrenormalised version.

{\bf Proof of Corollary \ref{lowdim}} To prove that for 
$1\leq p < 2$ (so $N=1$) we need make
no assumptions about the analytic continuation properties of
zeta functions, consider the expression for spectral flow obtained in 
Proposition \ref{killextraterms}. Specialising this formula to $N=1$ and 
so $k=0=j$ we find
\ben sf(\D,u^*\D u)C_{p/2+r}
 =-\frac{1}
{2}S\tau\left(q\{\tD,q\}(1+\tD^2)^{-(p/2+r)}\right)+holo,\een
where $holo$ is a function holomorphic at $r=(1-p)/2$. The super part of the 
trace gives
\be sf(\D,u^*\D u)C_{p/2+r}
=(-1/2)\tau( (u[\D,u^*]-u^*[\D,u])(1+\D^2)^{-p/2-r} )+holo\label{done}\ee
which, by the proof of Theorem \ref{grandfinale}, is (up to functions 
holomorphic at $r=(1-p)/2$) $(\sqrt{2\pi i})^{-1}\phi^r_1$ evaluated on 
$Ch_1(u^*)-Ch_1(u)$. Since $b\phi^r_1(a_0,a_1,a_2)$ is holomorphic at 
$r=(1-p)/2$ for all $a_0,a_1,a_2\in\A$, we can write Equation (\ref{done}) as
\ben sf(\D,u^*\D u)C_{p/2+r}=\frac{1}{\sqrt{2\pi i}}\phi^r_1(Ch(u))+holo.\een 
Taking residues gives
\ben sf(\D,u^*\D u)=\frac{1}{\sqrt{2\pi i}}\phi_1(Ch(u))\een
and the residue on the right necessarily exists and is given by $\phi_1$ by 
Proposition \ref{grandfinale}.

{\bf Remark}
In order to see how this last formula fits with Theorem 6.2 of \cite{CPS2},
we note that in \cite{CPS2} the spectral dimension is $1$ and that $p$ is a 
variable which we will write here as $p=1+z$. Then, (after switching the
roles of $u$ and $u^*$ and writing $\D$ in place of $D$) the formula in 
Theorem 6.2 of \cite{CPS2} is calculating:
$$1/2 res_{z=0}\tau(u^*[\D,u](1+\D^2)^{-1/2-z/2}).$$
While the formula in this paper is calculating:
$$res_{z=0}\tau(u^*[\D,u](1+\D^2)^{-1/2-z}).$$
It is a simple general fact about residues that these are identical.

\end{document}